\documentclass[a4paper]{elsarticle}

\usepackage[margin=3cm]{geometry}
\usepackage{mystyle}
\usepackage{multirow}
\usepackage{float}
\usepackage[table]{xcolor}

\usepackage[colorlinks,
            citecolor=blue,
            linkcolor=blue,
            urlcolor=black]{hyperref}

\journal{}

\newcommand{\K}{\bm{K}}
\renewcommand{\u}{\bm{d}} 
\newcommand{\f}{\bm{p}} 
\newcommand{\rhovec}{\bm{\rho}}
\newcommand{\x}{\bm{x}}

\definecolor{lightgray}{gray}{0.9}

\begin{document}

\title{Rank Reduction AutoEncoders for Mechanical Design: 
Advancing Novel and Efficient Data-Driven Topology Optimization}

\author[1]{Ismael Ben-Yelun\corref{cor1}}\ead{i.binsenser@upm.es}
\cortext[cor1]{Corresponding author}

\author[2]{Mohammed El Fallaki Idrissi}
\ead{mohammed.el_fallaki_idrissi@ensam.eu}

\author[3]{Jad Mounayer}
\ead{Jad.MOUNAYER@ensam.eu}

\author[4]{Sebastian Rodriguez}
\ead{sebastian.rodriguez_iturra@ensam.eu}

\author[2,3,4]{Francisco Chinesta}
\ead{Francisco.CHINESTA@ensam.eu}


\affiliation[1]{organization={E.T.S. de Ingeniería Aeronáutica y del Espacio, Universidad Politécnica de Madrid},
addressline={Pza. Cardenal Cisneros 3},%
city={Madrid},%
postcode={28040},%
country={Spain}%
}

\affiliation[2]{organization={ENSAM Institute of Technology
PIMM, ESI/Keysight Chair CNRS@CREATE, Singapore},
addressline={151 Bd de l'Hôpital},%
city={Paris},%
postcode={75013},%
country={France}%
}

\affiliation[3]{organization={ENSAM Institute of Technology
PIMM, SKF Chair},
addressline={151 Bd de l'Hôpital},%
city={Paris},%
postcode={75013},%
country={France}%
}

\affiliation[4]{organization={ENSAM Institute of Technology
PIMM, RTE Chair},
addressline={151 Bd de l'Hôpital},%
city={Paris},%
postcode={75013},%
country={France}%
}

\begin{abstract}
This work presents a data-driven framework for fast forward and inverse 
analysis in topology optimization (TO) by combining Rank Reduction 
Autoencoders (RRAEs) with neural latent-space mappings. The methodology 
targets the efficient approximation of the relationship between optimized 
geometries and their corresponding mechanical responses or Quantity of 
Interest (QoI), with a particular focus on compliance-minimized linear 
elastic structures. High-dimensional TO results are first compressed using 
RRAEs, which encode the data into a low-rank approximation via Singular 
Value Decomposition (SVD), obtained in this sense the most important 
features that approximate the data. Separate RRAE models are trained for 
geometry and for different types of QoIs, including scalar metrics, 
one-dimensional stress fields, and full two-dimensional von Mises stress 
distributions. The resulting low-dimensional latent coefficients of the 
latent space are then related through multilayer perceptrons to address both 
direct problems -- predicting structural responses from geometry -- and 
inverse problems---recovering geometries from prescribed performance 
targets. The proposed approach is demonstrated on a benchmark TO problem 
based on a half MBB beam, using datasets generated via density-based Solid 
Isotropic Material with Penalization (SIMP) 
optimization. Numerical results show that the framework enables accurate and 
computationally efficient surrogate models, with increasing robustness and 
fidelity as richer QoIs are considered. The methodology also provides a 
foundation for generative mechanical design by enabling the synthesis of new 
geometries and responses through latent-space exploration.
\end{abstract}


\maketitle


\section{Introduction}

Topology optimization (TO) has become a central tool in computational 
mechanics for the automated synthesis of high-performance structures under 
prescribed loading and boundary conditions. Since the seminal works on 
homogenization-based methods and density-based approaches, TO 
has been successfully applied to a wide range of engineering 
problems, including lightweight structural design 
\cite{li2018multi,junk2019structural,xue2021novel,zhang2021integrated}, 
compliant mechanisms 
\cite{zhang2018topology,sigmund1997design,zhu2020design}, and multi-physics 
systems \cite{park2012structural,munk2018multiobjective,das2020multi}. Among 
the available formulations, the Solid Isotropic Material with Penalization 
(SIMP) method \cite{bendsoe1989optimal} remains one of the most widely 
adopted techniques due to its conceptual simplicity and compatibility with 
Finite Element solvers \cite{reddy1993introduction,reddy2015introduction}.

Despite its success, TO remains computationally demanding, 
particularly when repeated analyses are required, such as in design space 
exploration, uncertainty quantification, real-time decision making, or inverse 
design. Each optimization run involves a large number of Finite Element 
problem resolutions, and the resulting optimized geometries are typically 
represented by high-dimensional density fields. This high computational cost 
has motivated the development of surrogate models and reduced-order 
representations capable of approximating optimized designs and their
mechanical responses at a fraction of the original cost.

Recent advances in machine learning, and deep learning in particular 
\cite{bengio2017deep}, have opened new opportunities for data-driven modeling 
in TO. Neural networks have been used to accelerate 
optimization, predict optimized layouts directly from problem parameters, and 
approximate solution fields such as displacements or stresses. Convolutional 
neural networks (CNNs) have been especially popular due to their ability to 
process grid-based density fields. However, purely black-box approaches often 
suffer from limited interpretability, poor extrapolation, and difficulties in 
handling inverse problems, where multiple geometries may correspond to the 
same performance metric.
An alternative and increasingly promising strategy consists of combining 
machine learning with reduced-order modeling techniques. Autoencoders, and 
more generally nonlinear manifold learning methods, provide a way to compress 
high-dimensional TO data into low-dimensional latent 
representations that capture the dominant geometric and physical features of 
the designs. These latent spaces can then be exploited for fast prediction, 
interpolation, and inverse analysis. Nevertheless, standard autoencoders 
often limit expressiveness or introduce redundancy and overfitting. 

In this context, Rank Reduction Autoencoders (RRAEs) 
\cite{mounayer2024rank} offer a flexible and physically motivated alternative 
by imposing that the latent space of an autoencoder should be expressed as 
low-rank approximation by means of a truncated Singular Value Decomposition 
(SVD) method \cite{klema1980singular,gerbrands1981relationships,mani1985application}. 
This architecture allows the model to identify the most relevant modes in the 
latent space while maintaining computational efficiency and robustness. The 
RRAE has demonstrated strong performance across a wide range of applications. 
Notable extensions include its variational autoencoder formulation, termed the 
VRRAE \cite{mounayer2025variational}. The RRAE has also been successfully 
applied to Generative Design (GD) and optimization of homogenized problems in 
composite materials \cite{idrissi2025new}, as well as to GD 
of parametric surrogate models \cite{idrissi2025generative}. More 
recently, the RRAE has been combined with DeepONets for heat transfer 
applications \cite{tierz2025variational}, employed in the modeling of complex 
dynamical systems \cite{mounayer2025rraedy} and the field of Structural 
Health Monitoring (SHM) for the detection of damage in plate structures using 
Lamb waves \cite{rodriguez2026damage}.

Therefore, the present work proposes a unified data-driven framework based on 
RRAEs to model the relationship between optimized geometries and their 
mechanical responses. Separate RRAE models are constructed for geometry and 
for different QoIs, including scalar stress measures, one-dimensional stress 
fields, and full two-dimensional von Mises stress distributions. The 
resulting low-dimensional latent coefficients are linked through neural 
regressors to address both direct problems—predicting responses from 
geometries—and inverse problems—reconstructing geometries from prescribed 
performance targets. It should be emphasized that the inverse problem is 
inherently ill-posed. Nevertheless, rather than applying a Newton-based 
approach to the forward problem, this work explores the construction of an 
inverse regression, fully acknowledging the risks and limitations associated 
with such a procedure. Through a detailed numerical study on a benchmark 
TO problem, the impact of the chosen QoI on model 
accuracy, robustness, and invertibility is systematically investigated.

The main contributions of this work are threefold:
(i) the introduction of RRAEs as an efficient and adaptive dimensionality 
reduction tool for TO data;
(ii) a comparative analysis of forward and inverse surrogate modeling 
performance as a function of the information content of the QoI; and
(iii) the demonstration of a latent-space-based framework that naturally 
enables GD and rapid exploration of optimized structural 
layouts.

The paper is structured as follows. Section \ref{sec:topology-optimization} 
provides an overview of the TO problem considered in this 
work for the generation of the dataset. Section \ref{sec:RRAE} introduces the 
main ingredients of the RRAE and how the architecture works. Following Section 
\ref{sec:methodology} presents how RRAE can be used to accelerate the topology 
optimization procedure and its use to perform inverse analysis. The 
performance of the proposed architecture in approximating as well as 
performing inverse analysis is demonstrated in Section \ref{sec:results} by 
considering a half MBB beam to perform TO and by 
considering as QoIs scalar, vector and 2d von Misses field. Finally, Section 
\ref{sec:concl_and_persp} provides conclusions and perspectives.

\section{Topology optimization}\label{sec:topology-optimization}

In this section, we review the topology optimization (TO) problem formulation of 
linear elastic solids using a Finite Element Method (FEM) 
\cite{reddy1993introduction,reddy2005introduction,reddy2015introduction} 
discretization for the case of compliance minimization problem with volume 
restriction.

\subsection{Problem formulation}

Starting in a rectangular domain $\Omega$ discretized by $n_e$ square finite elements and $n_d$ Degrees of Freedom (DoF), the equilibrium equation reads:
\begin{equation}
    \K(\rhovec) \u (\rhovec) = \f,
    \label{eq:equilibrium}
\end{equation}
where $\u\in \mathbb R^{n_d}$, $\f \in \mathbb R^{n_d}$, and 
$\K\in \mathbb R^{n_d \times n_d}$ are the displacements vector, the force 
vector, and the positive definite symmetric global stiffness matrix i.e., 
already reduced by the boundary conditions respectively 
\cite{reddy2015introduction}. The displacement vector and stiffness matrix 
depend on the element design variables, stored in the vector 
$\rhovec\in \mathbb R ^ {n_e}$, consisting of the relative densities of the 
elements. By defining the compliance $J$,
\begin{equation}
    J(\rhovec) := \f \cdot \u(\rhovec) = \u(\rhovec) \cdot \K(\rhovec) \u(\rhovec),
\end{equation}
as the objective function to be minimized, the optimization problem is given by:
\begin{subequations}
\begin{align}
\displaystyle\minimize_{\rhovec} & \quad J(\rhovec) = \u(\rhovec) \cdot \K(\rhovec) \u(\rhovec)= \sum_{e=1}^{n_e}\u_e(\rhovec) \cdot \K_e(\rhovec) \u_e(\rhovec)\label{eq:optimisation_subeq_cost}\\
 \mathrm{subject\;to} & \quad \vec K\left(\rhovec\right)\u \left(\rhovec\right) = \f \, , \label{eq:optimisation_subeq_equilibrium}\\
& \quad g(\rhovec) = \frac{V(\rhovec)}{V_{\max}} - f \leq 0 \, ,  \label{eq:optimisation_subeq_constraint} \\
& \quad \mathbf{0} \leq \rhovec \leq \vec I_{n_e} \,   . 
\end{align}
\label{eq:topology_optimization}
\end{subequations} 

Considering the area of the $e-$th element as $A_e$, the constraint 
$g(\rhovec)$ ensures that the volume $V(\rhovec) = \sum_e \rho_e A_e$ does not 
exceed a prescribed fraction $f$ of the maximum volume allowed in the domain 
$V_{\max} := V(\rhovec=\bm I_{n_e})$. The last constraint enforces the 
relative densities to acquire values between 0 and 1. Note that the force 
vector $\f$ is considered fixed throughout all the analyses. 

\subsection{Penalization and filtering}

In order to prevent numerical instabilities such as checker-boarding pattern 
\cite{sigmund1998numerical}, a stiffness penalization scheme is applied 
assuming Solid Isotropic Material with Penalization (SIMP) 
\cite{bendsoe1989optimal}. In this method, the Young's modulus $E_e$ of each 
element $e$ is related to its density $\rho_e$ through the following 
penalization function:
\begin{equation}
    E_e(\rho_e) = E_{\min} + \rho_e^p \left(E_0 - E_{\min} \right),
    \label{eq:penalization}
\end{equation}
where $E_0$ is the base material Young's modulus and $E_{\min}$ is a small 
numerical parameter to avoid the singularity of the stiffness matrix when 
$\rho_e \approx 0$. The penalization parameter $p$ is selected to enforce the 
elements to acquire either 0 or 1, and its value is typically set to $p = 3$. 
Due to the linear regime, the elements present a constant value of Young's 
modulus, and thus their local stiffness matrix $\K_e$ is linear on such value 
i.e., 
\begin{equation}    
    \K_e(\rho_e) = E_e(\rho_e)\K_0,
    \label{eq:stiffness_young}
\end{equation}
where $\K_0$ is the local stiffness matrix of a quad element without the 
influence of the Young's modulus and, hence, not dependant on the densities 
$\rhovec$. Then, the global stiffness matrix $\K$ is assembled. 

SIMP methods require filtering techniques, either density filter, sensitivity 
filter, or both. In order to obtain the filtered densities, a convolution with 
a cone (or, linear) kernel is applied to the (physical) densities $\rhovec$. 
We generate a weight matrix $w_{ij}$ by evaluating the centroidal coordinate 
pair $\bm c_i$ and $\bm c_j$ for the elements $i$ and $j$ as follows:
\begin{equation}
w_{ij} = \max\left\lbrace0, R- \lVert \vec c_i- \vec c_j \rVert_2 \right \rbrace \, , 
\end{equation}
where $R$ is the filter length prescribed and $||\cdot ||_2$ represents the 
Euclidean distance. With this filtering weights, the filtered density of an 
element $e$ is expressed as:
\begin{equation} \label{eq:filter}
\hat{\rho}_e = \dfrac{\sum_i w_{ei}\rho_i}{ \sum_i w_{ei}} \, .
\end{equation}

\subsection{Sensitivity analysis}

The TO problem set in \eqref{eq:topology_optimization} is 
solved using a gradient-based approach, e.g. the Method of Moving Asymptotes 
(MMA) \cite{svanberg1987method} or an Optimality Criteria (OC) method 
\cite{bendsoe1995optimization}. Therefore, the sensitivity analysis of the 
cost function $J(\rhovec)$ and the volume constraint $g(\rhovec)$ must be 
addressed. Note that the equilibrium equation constraint 
\eqref{eq:optimisation_subeq_equilibrium} is not considered in these 
approaches. Instead, the equilibrium is solved, and a local (convex) 
approximation around its solution is performed. Then, the convex optimization 
problem considering \eqref{eq:optimisation_subeq_cost} and 
\eqref{eq:optimisation_subeq_constraint} is solved at each iteration---see 
sequential approaches, e.g., \cite{christensen2008introduction}.

The gradient of the cost function \eqref{eq:optimisation_subeq_cost} is:
\begin{equation}
    \dfrac{\partial J(\rhovec)}{\partial \rho_e} =
    -\u(\rhovec) \cdot \dfrac{\partial \K(\rhovec)}{\partial \rho_e} \u(\rhovec),
    \label{eq:compliance_gradient}
\end{equation}
where it has been taken into account the derivative of the equilibrium equation \eqref{eq:optimisation_subeq_equilibrium} i.e.,
\begin{equation}
    \K(\rhovec) \dfrac{\partial \u(\rhovec)}{\partial \rho_e} = 
    - \dfrac{\partial \K(\rhovec)}{\partial \rho_e} \u(\rhovec).
\end{equation}

The gradient of the cost function \eqref{eq:compliance_gradient} can be expressed in a more convenient way as:
\begin{equation}
  \dfrac{\partial J(\rhovec)}{\partial \rho_e} =
  -p\left(E_0 - E_{\min}\right)^{p-1}\u_e(\rhovec)\cdot \K_0 \u_e(\rhovec), 
\end{equation}
where $\u_e$ is the displacement vector only considering the degrees of freedom of the element $e$, i.e., $\u_e \in \mathbb R^8$ in the 2D square element case. It has also been taken into account that:
\begin{equation}
    \dfrac{\partial \K_e(\rho_e)}{\partial \rho_e} = \dfrac{\partial E_e(\rho_e)}{\partial \rho_e}\K_0
    =
    p(E_0 - E_{\min})^{p-1} \K_0
\end{equation}
by computing the (straightforward) derivatives of equations \eqref{eq:penalization} and \eqref{eq:stiffness_young}. The gradient of the volume constraint \eqref{eq:optimisation_subeq_constraint} can be straightforwardly computed as:
\begin{equation}\textbf{}
    \dfrac{\partial g(\rhovec)}{\partial \rho_e} = A_e.
\end{equation}

In sensitivity filtering, a filter of sensitivities is computed using the chain rule, replacing the design variable $\rhovec$ with its filtered counterpart $\hat \rhovec$ i.e.,
\begin{equation}
    \dfrac{\partial J(\hat \rhovec)}{\partial \rho_e} = \sum_k \dfrac{\partial J(\hat \rhovec)}{\partial \hat \rho_k}\dfrac{\partial \hat \rho_k}{\partial \rho_e},\qquad
    \dfrac{\partial g(\hat \rhovec)}{\partial \rho_e} = \sum_k \dfrac{\partial g(\hat \rhovec)}{\partial \hat \rho_k}\dfrac{\partial \hat \rho_k}{\partial \rho_e}.
\end{equation}

Differentiating the filtering function \eqref{eq:filter} yields the term needed to compute the previous equation:
\begin{equation}
    \dfrac{\partial \hat\rho_k}{\partial \rho_e} = \dfrac{w_{ke}}{\sum_i w_{ki}}.
\end{equation}

\subsection{Stress computation}\label{eq:subsec-stress}

The von Mises stress distribution $\bm \sigma^{vM}\in \mathbb R^{n_e}$ is a 
vector containing the von Mises stress values $\sigma_e^{vM}$ at every 
element. This distribution is computed as a post-processing of the 
optimization problem, to account for the response (i.e., solution) of a 
structural problem in which the input is the geometry defined by the optimized 
topology.

First, the stress vector of an element $e$ as a function of the coordinate $\underline{\bm x}$ in the Voigt notation for the 2d case is:
\begin{equation*}
\bm \sigma_e(\underline{\bm x}) = \left(\sigma_{exx}, \sigma_{eyy}, \tau_{exy}\right)^T,
\end{equation*}
and it is computed in the FE analysis as follows:
\begin{equation}
    \bm \sigma_e(\underline{\bm x}) = \bm D_e \bm B_e(\underline{\bm x}) \bm d_e,
\end{equation}
where $\bm D_e$ is the plane stress constitutive matrix and $\bm B_e$ the 
kinematic matrix, both particularized for the element $e$. Note that the 
constitutive matrices may present different values since they depend on their 
element's stiffness $E_e(\rho_e)$. The von Mises stress measure reads:
\begin{equation}
    \sigma^{vM}_e(\underline{\bm x}) = \left(\sigma_{exx}^2+\sigma_{eyy}^2 -\sigma_{exx}\sigma_{eyy} +3\tau_{exy}^2 \right)^{\frac{1}{2}}.
\end{equation}

Lastly, to obtain a unique von Mises measurement per element $e$, $\sigma^{vM}_e(\underline{\bm x})$ is averaged in the element domain $\Omega_e$ performing numerical integration, namely:
\begin{equation}
    \sigma_e^{vM} = \int_{\Omega_e} \sigma_e^{vM}(\underline{\bm x})\mathrm{d}\Omega_e
    =
    \sum_{p=1}^{n_{ip}} \sigma_e^{vM}\left({\underline{\bm x}_p}\right)w_p^2,
\end{equation}
where $n_{ip}$ is the number of integration points in the Gaussian quadrature 
and $w_p$ the corresponding weights evaluated at such integration points $p$. 
Recall that quad 2d elements with $n_{ip}=4$ integration points are used in 
this paper.


\section{Rank Reduction Autoencoder (RRAE)}\label{sec:RRAE}

In this section, we review the Rank Reduction Autoencoder (RRAE) formulation, 
introduced by Mounayer et al. \cite{mounayer2024rank}, and highlight their 
advantages in comparison with regular autoencoders (AEs). This model is 
applied to perform dimensionality reduction to our datasets through a 
non-linear encoding. The input data can be of different nature, e.g., 
grid/pixel data, or a field defined in a domain, such as geometry or topology, 
among others. Let consider $\bm X = \lbrace \bm x_1, ..., \bm x_{n_s}\rbrace$, 
with $\bm x_j\in \mathbb R^T$ the input data consisting of $n_s$ samples. Note 
that the input space dimension can be either 1d or 2d i.e., 
$T = \left\lbrace D, D\times D\right\rbrace$, being $D$ the number of embedded 
features in the domain.

\begin{figure}[H]
    \centering
    \includegraphics[width=\linewidth]{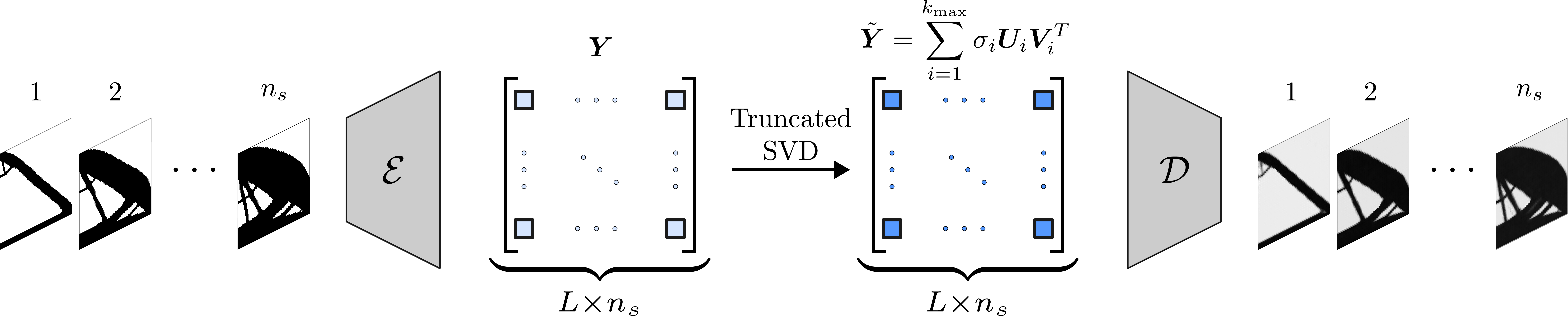}
    \caption{CNN RRAE architecture applied to a topology optimization 
    example, consisting of an encoder $\mathcal E$, performing the SVD and 
    decode $\mathcal D$. $L\gg T$ is the high-rank latent space dimension, 
    $k_{\max}$ is the low-rank latent size, i.e., number of modes retained 
    from the truncated SVD, and $n_s$ is the number of samples in the data.}
    \label{fig:rrae-geom}
\end{figure}

It is known that regular AEs perform the reconstruction through the following 
non-linear mapping $\tilde{\bm X} = \mathcal{D}(\mathcal{E}(\bm X))$ 
\cite{bengio2017deep}, where $\mathcal E(\cdot)$ and $\mathcal D(\cdot)$ are 
the encoder and decoder mappings, respectively. However, a low-rank 
approximation i.e., truncation, of a high-dimensional encoding mapping 
$\tilde{\bm Y}$ is introduced in RRAE formulation, to be decoded later, as 
shown in \figurename~\ref{fig:rrae-geom}.

Whereas AEs perform directly the map to a lower dimension, RRAEs perform the 
encoding mapping a sample $\bm x_j$ into a latent space of dimension $L$, and 
the concatenation of samples forms the matrix $\bm Y \in \mathbb R^{L\times n_s}$. 
Then, a low-rank approximation of this matrix is performed via a truncated 
Singular Value Decomposition (SVD), i.e., $\bm Y = \bm U \bm \Sigma \bm V^T$, 
where $\bm U \in \mathbb R^{L\times r}$, $\bm V^T \in \mathbb R^{r\times n_s}$, 
and the singular values are stored in the diagonal matrix 
$\bm \Sigma \in \mathbb R^{r\times r}$ sorted from highest to lowest such that 
$\bm \Sigma = \textrm{diag}(\sigma_i)$, $i=1,...,r$. Keeping only 
$k_{\max}\ll r$ modes from the singular values, so-called latent coefficients 
are defined and can be expressed as a linear combination of a orthogonal basis  
-- since $\bm Y$ is a real matrix -- consisting of the first $k_{\max}$ rows 
of the matrix $\bm U$. This basis is computed during the training process 
\cite{mounayer2024rank}.

Therefore, the reconstruction of a latent sample $\tilde{\bm y}_j$ after applying the truncated SVD might be expressed as:
\begin{equation}
    \tilde{\bm y}_j = \sum_{i=1}^{k_{\max}} \left(\sigma_i \bm U_i \bm V_i^T\right)_j =: 
    \sum_i^{k_{\max}} \alpha_{i,j} \bm U_i, 
    \quad
    j = 1,..., n_s,
\end{equation}
where $\alpha_{i,j} := \sigma_i \bm V_{i,j}^T$ are the above-mentioned latent coefficients. Thus, in matrix form, $\tilde{\bm Y}$ reads:
\begin{equation}
    \tilde{\bm Y} = \bm U \bm A,
\end{equation}
with $A_{i,j} = \alpha_{i, j}$, the matrix $\bm A \in \mathbb R^{k_{\max}\times n_s}$ 
containing these coefficients, and $\bm U \in \mathbb R^{L \times k_{\max}}$, 
with a slight abuse of notation, a truncated version of the $\bm U$ matrix 
from the SVD. Lastly, the reconstructed input data is performed applying the 
decoder map $\tilde{\bm X} = \mathcal D (\tilde{\bm Y})$. With that, a 
reconstruction loss accounting for the data-driven fitting is defined as 
follows:
\begin{equation}
    \mathcal L\left(\bm X,\tilde{\bm X}\right) = \dfrac{||\bm X-\tilde{\bm X}||_F}{|| \bm X||_F},
    \label{eq:reconstruction-loss}
\end{equation}
where $||\cdot ||_F$ corresponds to the Frobenius norm. An example of application of this architecture to 2d input data (images, thus CNN-RRAE) is depicted in \figurename~\ref{fig:rrae-geom}.

\section{Topology optimization using the RRAE}\label{sec:methodology}

This section presents the proposed data-driven framework for topology 
optimization based on Rank Reduction Autoencoders (RRAEs). The approach 
leverages the ability of RRAEs to compress high-dimensional topology 
optimization data into low-dimensional latent representations while adaptively 
identifying the most informative modes through rank truncation. By 
constructing separate reduced-order models for optimized geometries and for 
associated quantities of interest, the framework enables efficient learning of 
the relationships between design variables and structural responses. These 
latent representations are subsequently coupled through neural regressors to 
address both forward and inverse problems, providing a unified and 
computationally efficient surrogate modeling strategy for topology 
optimization.


Let $\bm X$ be a set of discretized geometry data, and each of the samples are 
introduced in a forward PDE problem e.g., a BVP such as a structural domain to 
solve the equilibrium, forming the set of discretized solutions $\bm S$ of 
such underlying physical problem. We will also refer to $\bm S$ as the 
Quantity of Interest (QoI) of the problem.

By applying an already trained RRAE, a geometry $\bm x_j$ can be expressed in the latent space from the couple $(\bm \alpha_j, \bm U)$ in the following way:
\begin{equation}
    \bm y_j = \mathcal E^g\left(\bm x_j\right) \quad \rightarrow \quad 
    \bm y_j = \bm U \bm \alpha_j \quad \Rightarrow \quad
    \bm\alpha_j = \bm U^T\bm y_j,
\end{equation}
where $\bm y_j$ is the associated embedded geometry and $\mathcal E^g(\cdot)$ is the encoder mapping of the geometries. The full set of samples might be expressed as $\bm A = \bm U^T\bm Y$, where $\bm Y = \bm U\bm A = \mathcal E^g(\bm X)$. Then, the reconstruction of a sample $\tilde{\bm x}_j$ is performed as:
\begin{equation}
    \tilde{\bm x}_j = \mathcal D^g\left(\bm U \bm \alpha_j\right),
\end{equation}
where $\mathcal D^g(\cdot)$ is the decoder mapping of the geometry RRAE.

Analogously, by separately training an RRAE for the QoI, a solution $\bm s_j$ can be expressed in its latent space given the so-defined couple $(\bm \beta_j, \bm V)$ of coefficient and basis such that:
\begin{equation}
    \bm z_j = \mathcal E^s(\bm s_j) \quad \rightarrow \quad
    \bm z_j = \bm V \bm \beta_j \quad \Rightarrow \quad
    \bm \beta_j = \bm V^T \bm z_j,
\end{equation}
where $\bm z_j$ is the embedded solution and $\mathcal E^s(\cdot)$ the encoder mapping for solutions. Similarly, the reconstruction of the solution $\bm s_j$ is:
\begin{equation}
    \tilde{\bm s}_j = \mathcal D^s\left(\bm V\bm \beta_j\right),
\end{equation}
where $\mathcal D^s(\cdot)$ is the decoder mapping for the solutions.

In other words, low-rank approximations for the geometries in terms of latent coefficients $\bm \alpha_j$ and basis $\bm U$, and for the solutions i.e., QoIs, in terms of latent coefficients $\bm \beta_j$ and basis $\bm V$ are found. In order to join both latent spaces, mapping functions relating $\bm\alpha_j$ to $\bm \beta_j$ and vice versa have to be obtained. Therefore, MLPs are used as interpolators to build these functions. Specifically, we will refer to $\mathcal{NN}^d_{\star}: \bm \alpha_j \rightarrow \bm \beta_j$ as the neural network addressing the `direct' problem i.e., from geometry to solution, and $\mathcal{NN}^i_{\star}: \bm \beta_j \rightarrow \bm \alpha_j$ as the neural network connecting the `inverse' problem i.e., from solution to geometry.

Considering the joint of these two architectures, two novel frameworks are developed: fast solution inference for a (potentially new) geometry -- so-called `direct' problem -- and geometry discovery from a target performance---so-called `inverse' problem. In the first case, the steps to apply the direct problem given a geometry $\bm x_j$ are:
\begin{itemize}
    \item Geometry embedding: $\bm y_j = \mathcal E^g(\bm x_j )$
    \item Projection in the (geometry) latent space: $\bm \alpha_j = \bm U^T\bm y_j$
    \item Evaluate the associated coefficients $\bm \beta_j$ using the trained regression: $\bm \beta_j = \mathcal{NN}^d_{\star}(\bm \alpha_j)$
    \item Solution evaluation: $\tilde{\bm s}_j = \mathcal D^s(\bm V \bm \beta_j)$
\end{itemize}

In the second case, the steps to address the inverse problem are analogous. Starting from a solution $\bm s_j$,
\begin{itemize}
    \item Solution embedding: $\bm z_j = \mathcal E^s(\bm s_j)$
    \item Projection in the (solution) latent space: $\bm \beta_j = \bm V^T \bm z_j$
    \item Evaluate the latent coefficient using (another) trained regression:  $\bm\alpha_j = \mathcal{NN}^i_{\star}(\bm\beta_j)$
    \item Geometry retrieval: $\tilde{\bm x}_j = \mathcal D^g(\bm U \bm\alpha_j)$
\end{itemize}

We refer to these two ways as the following mappings:

\begin{equation*}
\tilde{\bm S} = \texttt{geo2sol}_{\star}(\bm X),\qquad
\tilde{\bm X} = \texttt{sol2geo}_{\star}(\bm S),
\end{equation*}

\begin{figure}[H]
    \centering
    \subfloat[Scalar QoI.\label{fig:rrae-geom-and-scalar}]
    {\includegraphics[width=.8\linewidth]{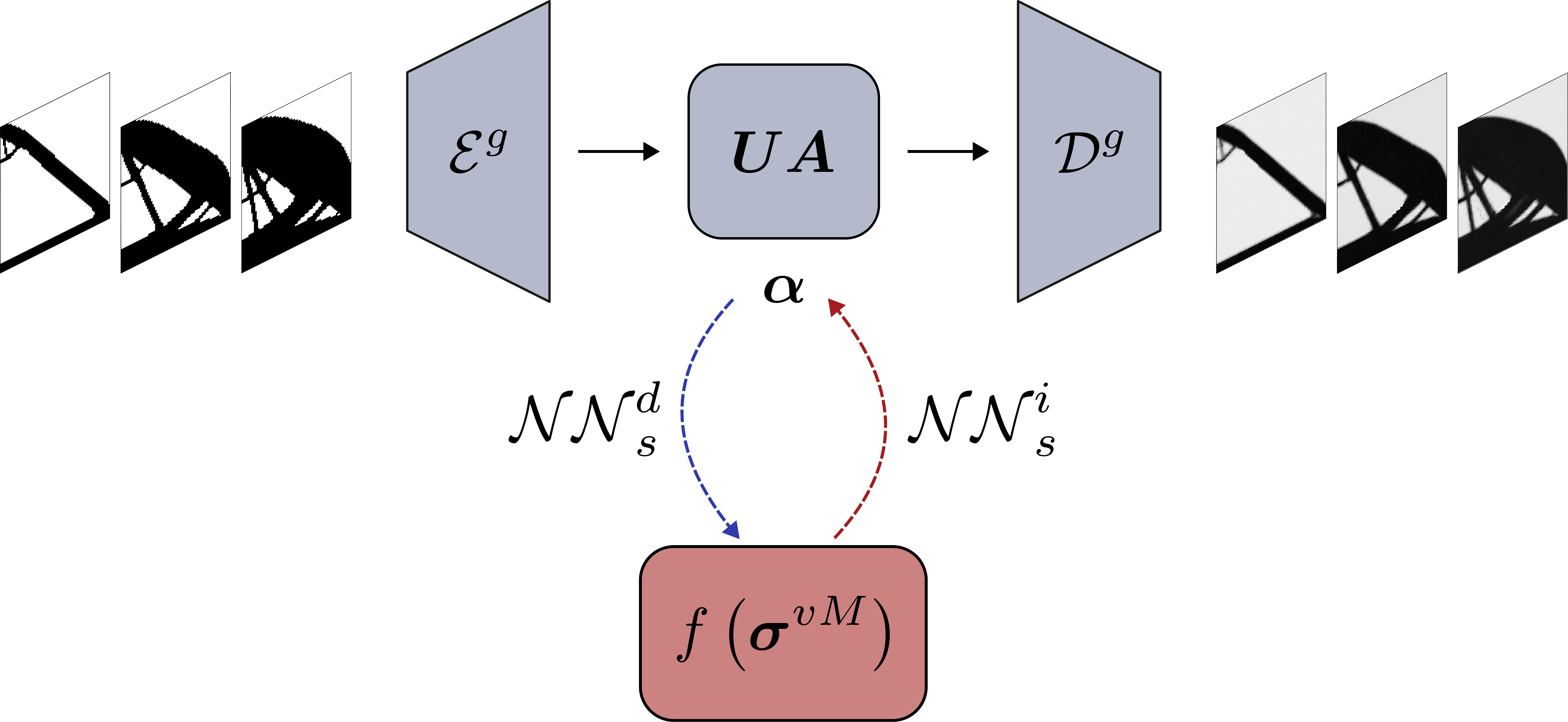}} \\
    \subfloat[2d field in a fixed domain $\Omega$ QoI.\label{fig:rrae-geom-and-sol2d}]{\includegraphics[width=.8\linewidth]{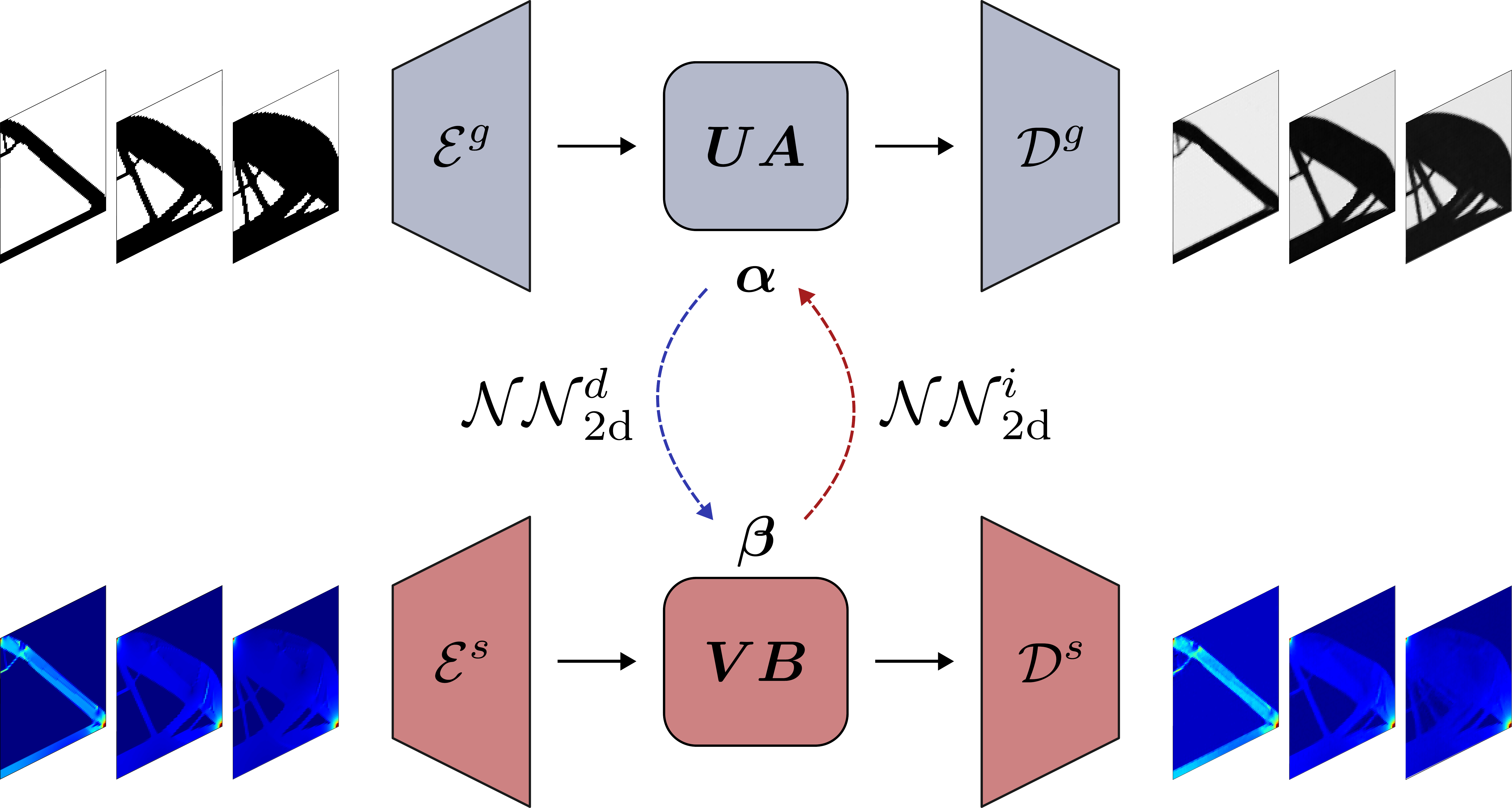}}
    \caption{Pipeline to address the direct \texttt{geo2sol}$_{\star}$ or 
    inverse \texttt{sol2geo}$_{\star}$ problems, for (a) scalar and (b) 2d 
    field QoIs. There are three models in each case that are trained 
    \textit{separately}: geometry RRAE $(\bm U\bm A)$, solution RRAE 
    $(\bm V \bm B)$, and the mappings between both latent coefficients, either 
    $\mathcal{NN}^d_{\star}:\bm\alpha\rightarrow \bm\beta$ or 
    $\mathcal{NN}^i_{\star}: \bm\beta \rightarrow \bm\alpha$.}
    \label{fig:rrae-geom-and-sol}
\end{figure}

While keeping the geometry RRAE fixed, in this paper we propose three 
different QoIs regarding the solutions, namely:
\begin{itemize}
    \item \textbf{Scalar}. This QoI does not require a low-rank approximation, 
    therefore the mapping directly relates the (geometry) latent coefficients 
    with the scalar value. In this case, the scalar is a function of the von 
    Mises stress distribution i.e., $f(\bm \sigma^{vM})$, as highlighted in 
    \figurename~\ref{fig:rrae-geom-and-scalar}.
    \item \textbf{1d field in a fixed domain}. This QoI comprises the von 
    Mises stress values of the elements belonging in the diagonal of the 
    domain $\Omega$, since this is the main load path. Since a curve is 
    obtained as the QoI -- hence the 1d field -- an MLP-RRAE is used as the 
    low-rank model. 
    \item \textbf{2d field in a fixed domain}. The QoI is the whole von Mises 
    distribution in the domain, $\bm \sigma^{vM}$ expressed in a matrix form 
    to be understood as an image by the CNN-RRAE. This is depicted in 
    \figurename~\ref{fig:rrae-geom-and-sol2d}. 
\end{itemize}

Lastly, another advantage that this methodology might present is its use in 
GD. New geometries or solutions might be generated by 
choosing new points in the latent space, i.e., $\hat{\bm\alpha}$ or 
$\hat{\bm \beta}$ and then decoding in the corresponding problem---geometry or 
solution, respectively. Particularly, these new points might be expressed as 
interpolations of the training set latent coefficients, i.e., within the 
convex hull of the latent coefficients space in order to avoid extrapolations. 
Although this is a potentially powerful tool for GD with applications such as 
dataset enriching, discovery of new topologies, etc., this proposal is out of 
the scope of this paper. 

\section{Numerical results}\label{sec:results}

\begin{figure}[H]
    \centering
    \includegraphics[width=0.25\linewidth]{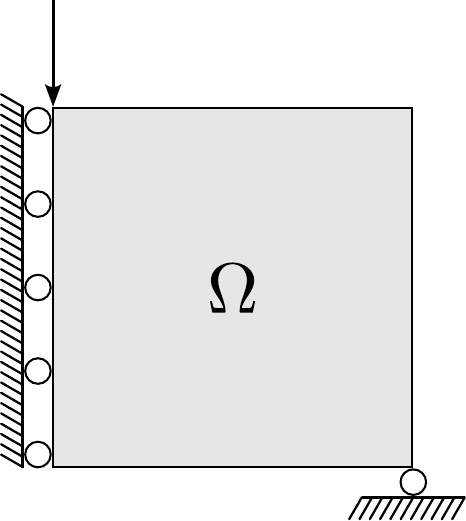}
    \caption{Domain $\Omega$ and boundary conditions of the half MBB beam for the TO problem.}
    \label{fig:boundary-conditions}
\end{figure}

We now proceed to demonstrate the utility and efficiency of the proposed 
methodology using a TO example and different QoIs. The main example consists 
of performing a Design of Experiments (DoE) by sweeping the volume constraint 
in a compliance minimization problem, as the half MBB beam displayed in 
\figurename~\ref{fig:boundary-conditions}.

The domain $\Omega$ is discretized into $80\times 80$ quad FEs of linear 
elastic base material with unit Young's modulus $E_0 =1$ and Poisson's ratio 
$\nu = 0.3$. For the TO problem, a SIMP scheme with a penalization parameter 
of $p=3$ and a sensitivity filter as the one described in 
Section~\ref{sec:topology-optimization}, with a filter radius $R = 1.5$ of the 
element's size is applied. A minimum Young's modulus value of 
$E_{\min} = 10^{-9}$ is introduced in the elements  in order to avoid 
obtaining singular stiffness matrices $\bm K$. Regarding the volume constraint 
$f$, this value has been varied by generating a vector of $N=100$ elements 
equispaced linearly between 0.1 and 0.9. Therefore, a total of $N$ TO problems 
are solved using the Optimality Criteria (OC) method provided in the 
\texttt{top88.m} open source script by DTU \cite{andreassen2011efficient} in 
order to generate the dataset to be used in this proposed methodology. Note 
that the boundary conditions (BC) are fixed in the analysis, wherein the force 
applied is a unit load. There is a possibility to parameterize the BCs e.g., 
by displacing the load application and/or support node along the boundary of 
the domain $\partial\Omega$, however, this idea is out of the scope of this 
paper.

Once the set of optimal topologies $\lbrace\rhovec^*_1,...,\rhovec_N^*\rbrace$ 
is obtained, they are post-processed and the corresponding set of von Mises 
stress distributions $\lbrace\bm \sigma^{vM,*}_1,...,\bm \sigma^{vM,*}_N\rbrace$ 
are computed according to Section~\ref{eq:subsec-stress}. The set of optimal 
topologies constitutes the discretized geometries dataset $\bm X$, in which a 
sample $\bm x_j$ is defined by:
\begin{equation*}
\bm x_j := \rhovec_j^*,
\end{equation*}
properly rearranged in a matrix i.e., replicating the 2d mesh. On the other 
hand, the von Mises stress distribution set constitutes the discretized 
solutions dataset $\bm S$, and each sample $\bm s_j$ is assigned as a function 
of the $\bm \sigma^{vM,*}_j$ distribution depending on the type of QoI 
considered $\lbrace$scalar, 1d, 2d$\rbrace$.

Thus, the training data $\bm X$ and $\bm S$ are generated to train the 
geometries and solutions RRAEs separately. Afterwards, we train models that 
represent a mutual fitting between both latent spaces. In addition to these 
datasets, a test dataset is generated to assess the accuracy of the 
methodology in new data samples. To this end, another DoE of $N_{\rm test} = 20$ 
volume fractions linearly equispaced between 0.1 and 0.9 is performed. Note 
that there are no common volume fraction values with respect to the training 
data but 0.1 and 0.9.

In the following subsections, three types of QoIs are analyzed: scalar, 1d 
field in a fixed domain, and 2d field in a fixed domain. In each of them, the 
direct and inverse problems are addressed.

\subsection{Scalar quantity of interest}\label{subsec:scalar-qoi}

The first case corresponds to the relation between the geometries latent 
coefficients and an equivalent or global response of the structure. This 
global magnitude is the maximum von Mises stress value in the domain, which 
may be of particular interest for industrial applications since it is often 
included in TO problems as a constraint. Thus, the samples $s_j$ of the 
solutions dataset $\bm S$ are assigned as follows:
\begin{equation*}
s_j := \max\left(\bm\sigma^{vM,*}_j\right).
\end{equation*}

These scalar problems are the most challenging of the three. The more 
compressed the solution is -- a global magnitude contains less information 
than the full resolution field -- the less bijective the problem is. Hence, 
the inverse problem becomes more ill-posed.

The next step is to build the geometries RRAE model. As the input are the 
discretized relative densities $\rhovec_j^*$ of the elements, rearranged in 
matrix form to maintain the mesh layout and neighborhood i.e., an image, a 
convolutional autoencoder architecture is used (CNN-RRAE). The chosen latent 
size is $L=500$, to further compress via SVD retaining $k_{\max} = 2$ modes. 
Theoretically, an accurate reconstruction can be performed with only one 
mode -- since it is only the volume fraction the parameter varying across the 
dataset -- however, we keep 2 modes in order to obtain predictions with higher 
accuracy. The rest of hyper-parameters are displayed in 
\tablename~\ref{tab:geometries-rrae}.
\begin{table}[H]
    \centering
    \begin{tabular}{c|ccccc}
    \hline
    Module & Layer & Input shape & Output shape & Activation & Kernel \\ \hline
    \rowcolor{lightgray}\cellcolor{white} & Conv. Layer & $(n_b, 1, 80, 80)$ & $(n_b, 32, 40, 40)$ & ReLU & $(3\times 3)$ \\ 
    \rowcolor{lightgray}\cellcolor{white}& Conv. Layer & $(n_b, 32, 40, 40)$ & $(n_b, 64, 20, 20)$ & ReLU & $(3\times 3)$\\
    \rowcolor{lightgray}\cellcolor{white}& Conv. Layer & $(n_b, 64, 20, 20)$ & $(n_b, 128, 10, 10)$ & ReLU & $(3\times 3)$ \\
    \rowcolor{lightgray}\cellcolor{white}\multirow{-4}{*}{Encoder} & Dense & $(n_b, 12\,800)$ & $(n_b, 500)$ & ReLU & $-$ \\
    & Dense & $(n_b, 500)$ & $(n_b, 3200)$ & ReLU & $-$\\
    & Conv. Trans. Layer & $(n_b, 32, 10, 10)$ & $(n_b, 128, 20, 20)$ & ReLU & $(3\times 3)$ \\
    Decoder& Conv. Trans. Layer & $(n_b, 128, 20, 20)$ & $(n_b, 64, 40, 40)$ & ReLU & $(3\times 3)$\\
    & Conv. Trans. Layer & $(n_b, 64, 40, 40)$ & $(n_b, 32, 80, 80)$ & ReLU & $(3\times 3)$\\
    & Final Conv. Layer & $(n_b, 32, 80, 80)$ & $(n_b, 1, 80, 80)$ & Linear & $(1\times 1)$\\
    \hline
    \end{tabular}
    \caption{Hyper-parameters of the CNN-RRAE model for the geometries, where 
    $n_b$ is the number of samples in a batch. The parameters common to all 
    convolutional layers in both encoder and decoder are 
    ${\rm padding} = 1$, ${\rm stride} = 2$ and ${\rm dilation} = 1$, except 
    for the final convolutional layer, which has no padding and 
    ${\rm stride} = 1$ in order to recover the initial $80\times 80$ 
    shape---acting as a `trainable' pooling for the channels.}
    \label{tab:geometries-rrae}
\end{table}

Since it is already defined in $[0, 1]$, geometry data $\bm x_j$ is fed to the 
CNN-RRAE model without any normalization. The training parameters are shown in 
\tablename~\ref{tab:geometries-rrae-optimization-parameters}. There are three 
training stages, in which the learning rate is decreased to a tenth of the 
value of the previous stage so smoother steps are performed at later epochs.
\begin{table}[H]
    \centering
    \begin{tabular}{c|c}
         Parameter & Choice \\ \hline
         Optimizer & AdaBelief \\
         \rowcolor{lightgray}Learning rate & $[10^{-3}, 10^{-4}, 10^{-5}]$\\
         Epochs & $[3500, 3500, 3500]$ \\
         \rowcolor{lightgray}Batch size & $[20, 20, 20]$ \\
         Loss & Equation~\eqref{eq:reconstruction-loss}
    \end{tabular}
    \caption{CNN-RRAE geometries model. Optimization parameters of the training.}
    \label{tab:geometries-rrae-optimization-parameters}
\end{table}

\begin{figure}[H]
    \subfloat[Train]
    {\includegraphics[width=.45\linewidth]{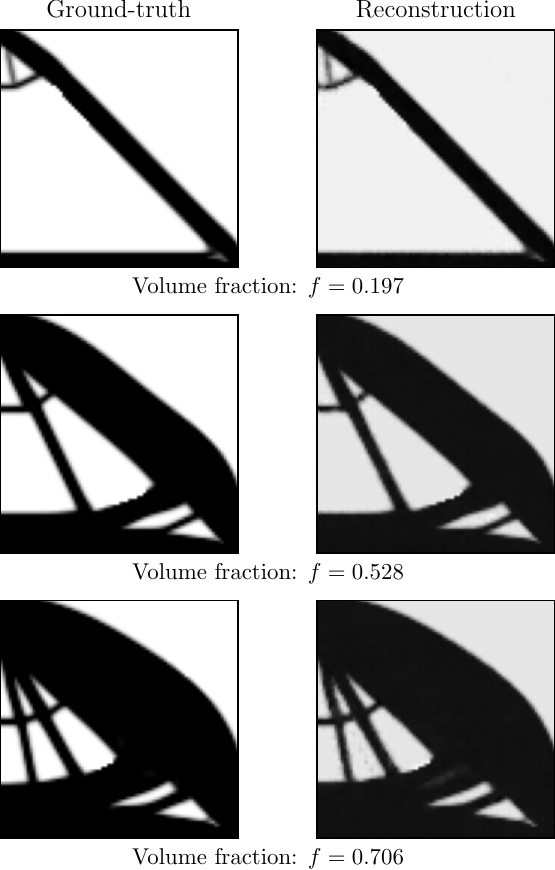}}\hfill%
    \subfloat[Test]
    {\includegraphics[width=.45\linewidth]{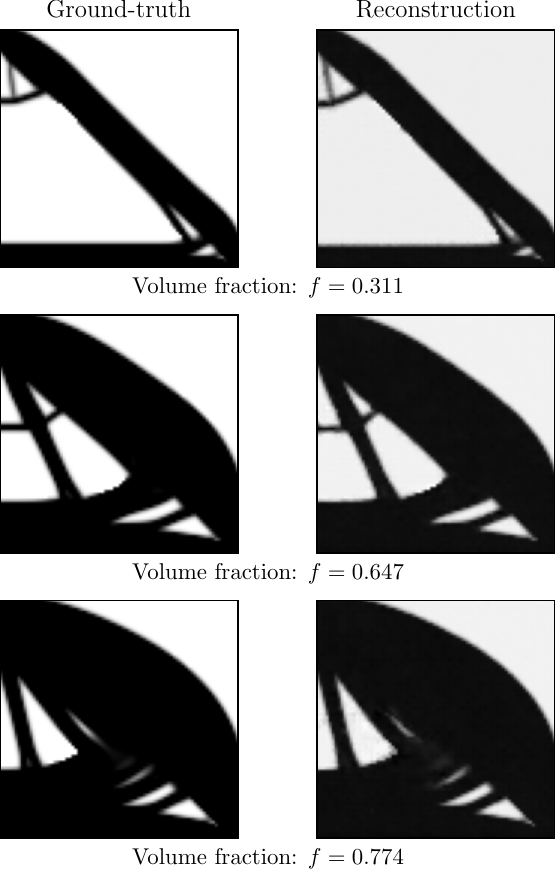}}
    \caption{CNN-RRAE geometries model. Example reconstructions in (a) train and (b) test sets.}
    \label{fig:rrae-geom-reconstruction}
\end{figure}

After the training, the reconstruction loss of the geometries RRAE is 
$\mathcal L_{\rm train} = 4.32\,\%$ in the training set, and its value in the 
test counterpart is $\mathcal L_{\rm test} = 6.32\,\%$. Examples for both 
train and test reconstructions for three arbitrary samples are depicted in 
\figurename~\ref{fig:rrae-geom-reconstruction}.

\subsubsection{Direct problem}

Once obtained the geometries latent coefficients $\bm A$, we propose a mapping 
relating such coefficients with the scalar QoI i.e., the maximum von Mises 
stress $\sigma_{\max}^{vM}$ in the domain. In order to build the whole 
pipeline within a DL framework, an MLP is trained to this end. Thus, the 
sought mapping can be expressed as 
$\mathcal{NN}_s^d:\bm \alpha_j\rightarrow  s_j$, where the superscript $d$ 
stands for direct problem. In order to only highlight the hyper-parameters of 
the RRAE models -- those of interest in this paper -- the hyper-parameters of 
the MLP architecture are displayed in \tablename~\ref{tab:nnd-1d} in 
\ref{app:nnd-1d}. Both inputs and outputs are normalized with an standard 
scaler i.e., each feature being subtracted its mean and divided by its 
standard deviation. An 80-20\% split is performed to generate the validation 
set. Then, the training is carried out with the optimization parameters 
highlighted in \tablename~\ref{tab:nnd-1d-optimization-parameters}, also in 
\ref{app:nnd-1d}.

Performing the training, the loss value in the training set is 
$\mathrm{MSE}_{\rm train} = 1.4\cdot 10^{-1}$, which represents a slightly 
high value for this kind of standards. This may be also observed in the 
coefficient of determination, obtaining $R^2_{\rm train} = 0.814$ and 
$R^2_{\rm test} = 0.930$. The fitting curve comparing the ground-truth 
$\sigma_{\max}^{vM}$ against the $\mathcal{NN}_s^{d}$ predicted values is 
depicted in \figurename~\ref{fig:nnd-s-latent}, in \ref{app:mlp-fitting-curves}. 
All the fitting curves from the $\mathcal{NN}^d_{\star}$ and 
$\mathcal{NN}^i_{\star}$ MLP models are depicted in such appendix to avoid the 
proliferation of plots.

The next step within the scalar QoI direct problem is to join the CNN-RRAE 
geometries and MLP $\mathcal{NN}_s^{d}$ models to create a pipeline. Such 
pipeline is referred to as $\texttt{geo2sol}_s$, performing the following 
steps (also described in Section~\ref{sec:methodology}): (1) takes a geometry 
$\bm x_j$ as input, (2) encodes such geometry through the first part of the 
CNN-RRAE geometry model to obtain the corresponding latent coefficient 
$\bm\alpha_j$, and (3) maps such latent coefficient through 
$\mathcal{NN}_s^d$ to predict the scalar QoI $\tilde{s}_j$, i.e., the maximum 
von Mises stress $\sigma_{\max}^{vM}$ in the original input geometry. Thus, 
following this procedure allows solutions to be obtained without having to go 
through the FE solver, and therefore achieving a surrogate model enabling 
computational efficiency.

To evaluate the accuracy of predictions in train and test datasets, the 
coefficient of determination between the ground-truth $\sigma_{\max}^{vM}$ 
values and the $\texttt{geo2sol}_s$ predicted ones is computed, yielding the 
following:
\begin{equation*}
    R^2_{\rm train} = 0.867, \quad R^2_{\rm test} = 0.987.
\end{equation*}

\begin{figure}[H]
    \centering
    \subfloat[Train]{\includegraphics[width=0.475\linewidth]{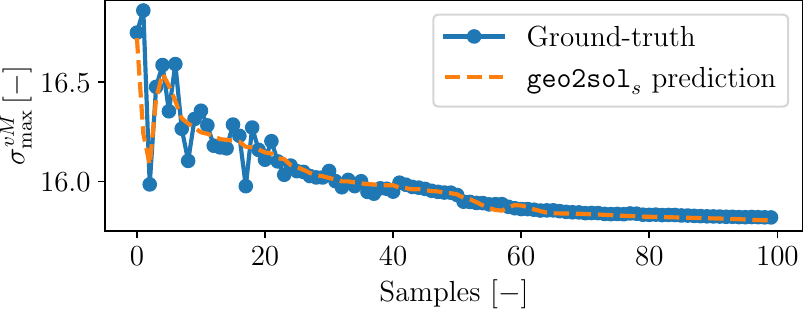}}\hfill
    \subfloat[Test]{\includegraphics[width=0.475\linewidth]{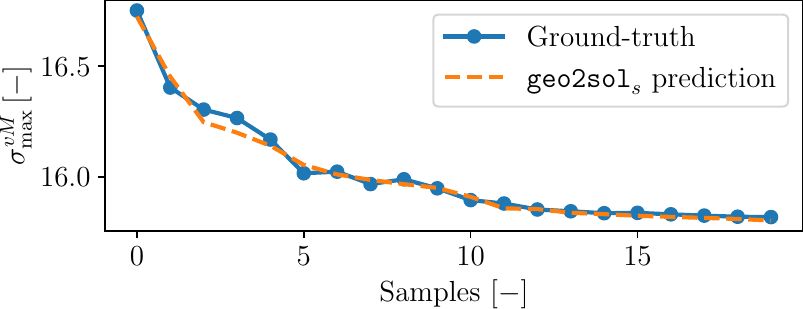}}
    \caption{$\texttt{geo2sol}_{s}$ predictions in (a) train and (b) test.}
    \label{fig:geo2sol-s}
\end{figure}

Since the predictions in the $\mathcal{NN}_s^d$ are not sufficiently accurate, 
this is reflected in the $\texttt{geo2sol}_s$ predictions, especially for the 
train set. As it was mentioned, this case presents worse performance due to 
non-bijectivity of the $\texttt{geo2sol}_s$ mapping. The train and test curves 
displaying the true values against the predicted ones with this direct 
methodology for the scalar QoI are depicted in \figurename~\ref{fig:geo2sol-s}.

\subsubsection{Inverse problem}

In this subsection, an analogous inverse procedure is performed. In this case, 
we propose a way of obtaining geometries from the scalar QoI. In order to do 
so, the first step is to build the inverse mapping 
$\mathcal{NN}_s^i:s_j\rightarrow \bm \alpha_j$, where the superscript $i$ 
stands for inverse. This model is an MLP with the same architecture and 
training parameters as its direct counterpart -- Tables \ref{tab:nnd-1d} and 
\ref{tab:nnd-1d-optimization-parameters} -- with the input and output 
dimensions permuted since $s_j\in \mathbb R^1$ and 
$\bm \alpha_j\in\mathbb R^2$. Normalizing both inputs and outputs, and 
training the model, a loss value of $\mathrm{MSE}_{\rm train} = 6.4\cdot 10^{-2}$ 
is obtained. The coefficients of determination between ground-truth and 
$\mathcal{NN}_s^i$ predicted latent coefficients values 
$\bm \alpha = [\alpha_0, \alpha_1]$ are $R^2_{\rm train} = 0.943$ and 
$R^2_{\rm test} = 0.945$. The corresponding curves are depicted in 
\figurename~\ref{fig:nni-s-latent} (\ref{app:mlp-fitting-curves}).

Having fitted $\mathcal{NN}_s^i$, the pipeline $\texttt{sol2geo}_s$ (inverse 
to the previous pipeline) can be generated. In this case, a certain value of 
$\sigma_{\max}^{vM}$ i.e., the scalar QoI $s_j$, is requested, so it is fed to 
the pipeline as an input. The value $s_j$ is mapped into the corresponding 
geometry latent coefficient $\bm \alpha_j$ through $\mathcal{NN}_s^i$. 
Finally, a geometry $\tilde{\bm x}_j$ is retrieved by passing $\bm \alpha_j$ 
through the decoder of the CNN-RRAE geometries model. Since the output is a 
mesh that can be interpreted as an image with each pixel having a density 
value $\bm \rho_j^*\in\left[0,1\right]$, the reconstruction error between the 
ground-truth geometry $\bm x_j$ and the $\texttt{sol2geo}_s$ predicted 
geometry $\tilde{\bm x}_j$ is performed. This in turn implies the computation 
of the $L^2$ norm of the difference of the flattened matrices, as highlighted 
in Equation~\eqref{eq:reconstruction-loss}, yielding the following loss values:
\begin{equation*}
\mathcal{L}_{\rm train} = 20.59\,\%, \quad \mathcal{L}_{\rm test} = 20.22\,\%.
\end{equation*}
\begin{figure}[H]
    \centering
    \subfloat[Train example]{\includegraphics[width=.8\linewidth]{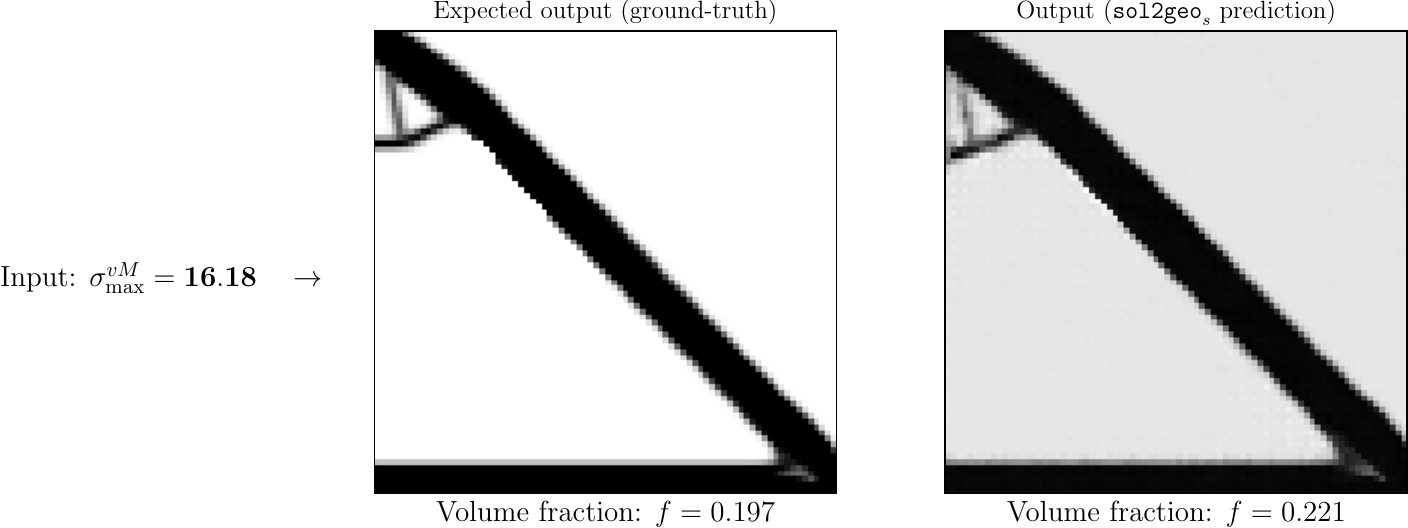}}\\
    \subfloat[Test example]{\includegraphics[width=.8\linewidth]{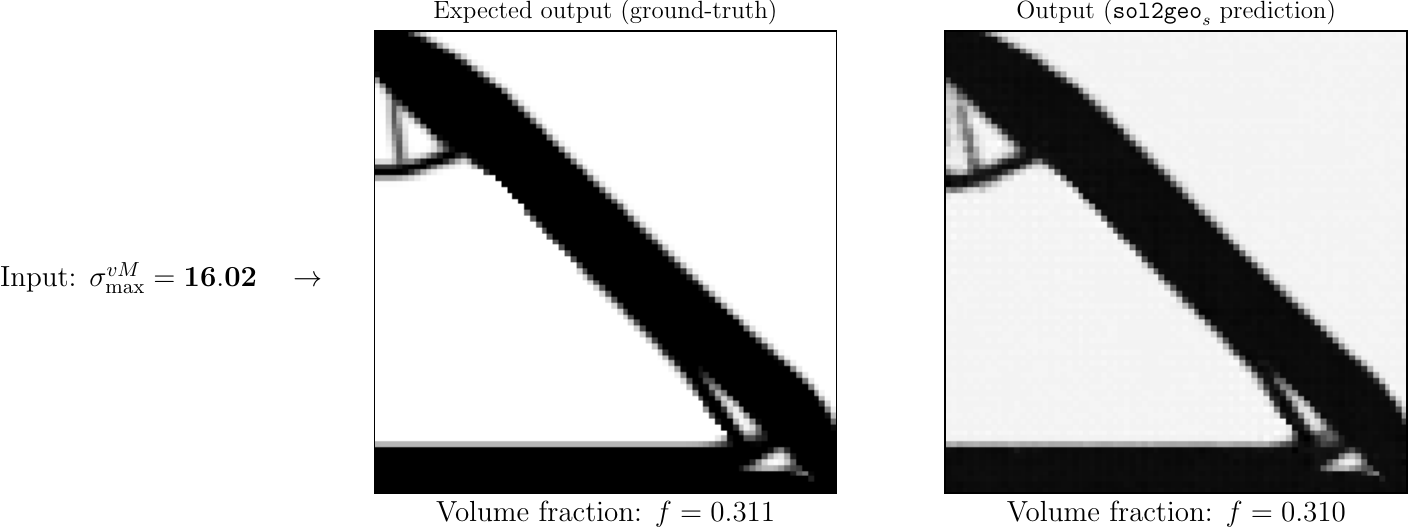}}
    \caption{$\texttt{sol2geo}_s$ example predictions in (a) train and (b) test.}
    \label{fig:sol2geo-s}
\end{figure}

\figurename~\ref{fig:sol2geo-s} displays examples of $\texttt{sol2geo}_s$ 
train and test predictions when requesting a $\sigma_{\max}^{vM}$ value, 
against the corresponding ground-truth geometries. As can be noticed in the 
errors or in the figure, it is important to emphasize that actually in this 
scalar case architecture the model is learning in an overfitted way, i.e., 
memorizing, a `bijective' type of mapping, which is not a good practice for 
this type of models. This highlights that more information is needed on the 
solution side for this methodology to be successful, as we shall show in the 
next two cases in this paper.  Additionally, when the requested value is out 
of the bounds of the training set, this pipeline starts to produce geometries 
that lack physical meaning, which is another drawback in this first case.

\subsection{1d field quantity of interest in a fixed domain}\label{subsec:1dfield-qoi}

This second analysis showcases the relation between the geometries latent 
coefficients (from the previous model) and the latent coefficients of the 
response of the structure in a 1d field. In this case, the von Mises stress 
values of the elements of the diagonal (top-left to bottom-right elements in 
$\Omega$). They are of special interest since they comprise the main load 
path---connecting to nodes affected by the BC. Thus, the samples $\bm s_j$ of 
the solutions dataset $\bm S$ are assigned as follows:
\begin{equation*}
\bm s_j := \mathrm{diag}\left(\bm \sigma^{vM,*}_j\right).
\end{equation*}

The RRAE model regarding the geometries is the same as the previous 
Section~\ref{subsec:scalar-qoi}. The next step is to build the solutions RRAE 
model. The solutions input is dense i.e., a matrix, therefore an MLP-RRAE 
architecture is selected. The latent size is $L = 500$, and the truncated SVD 
compression is performed keeping $k_{\max} = 1$. The rest of hyper-parameters 
are displayed in \tablename~\ref{tab:solutions-rrae}.
\begin{table}[H]
    \centering
    \begin{tabular}{c|cccc}
    \hline
    Module & Layer & Input shape & Output shape & Activation \\ \hline
    \rowcolor{lightgray}\cellcolor{white} & Dense & $(n_b, 80)$ & $(n_b, 64)$ & ReLU \\ 
    \rowcolor{lightgray}\cellcolor{white}\multirow{-2}{*}{Encoder} & Dense & $(n_b, 64)$ & $(n_b, 500)$ & Linear \\
    & Dense \#1 & $(n_b, 500)$ & $(n_b, 64)$ & ReLU\\
    & Dense \#2 & $(n_b, 64)$ & $(n_b, 64)$ & ReLU\\
    Decoder & & $\vdots$ & & \\
    & Dense \#6& $(n_b, 64)$ & $(n_b, 64)$ & ReLU\\
    & Final Layer & $(n_b, 64)$ & $(n_b, 80)$ & Linear\\
    \hline
    \end{tabular}
    \caption{Hyper-parameters of the MLP-RRAE model for the von Mises stress 
    1d field solutions, where $n_b$ is the number of samples in a batch. In 
    the decoder, the properties of dense layers \#2 to \#6 are the same, hence 
    they are omitted for the sake of simplicity.}
    \label{tab:solutions-rrae}
\end{table}

\begin{figure}[H]
    \subfloat[Train]
    {\includegraphics[width=.4\linewidth]{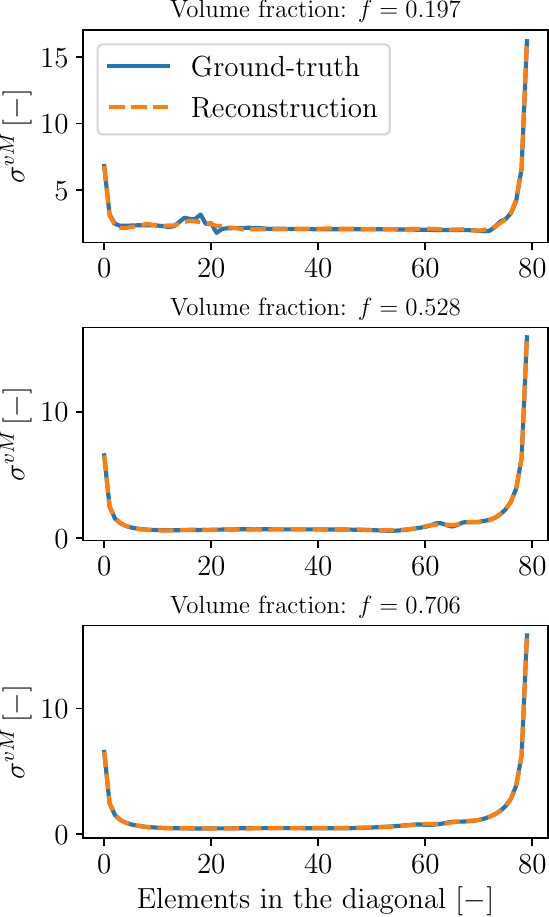}}\hfill%
    \subfloat[Test]
    {\includegraphics[width=.4\linewidth]{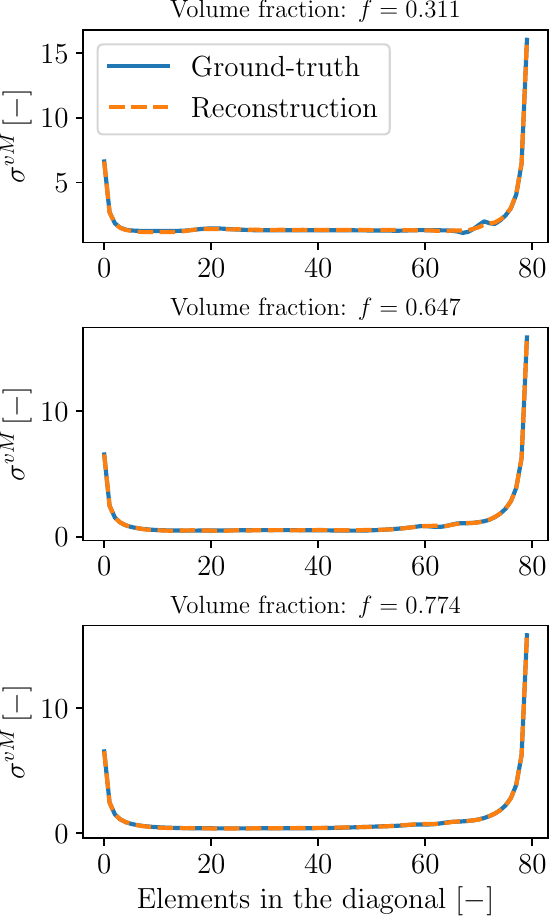}}
    \caption{MLP-RRAE 1d solutions model (von Mises stress in the diagonal). Example reconstructions in (a) train and (b) test sets.}
    \label{fig:rrae-sol-reconstruction}
\end{figure}

The training parameters are the same as those displayed in 
\tablename~\ref{tab:geometries-rrae-optimization-parameters}. Once fitted, the 
reconstruction loss of the solutions RRAE is $\mathcal L_{\rm train} = 3.14\,\%$ 
in the training set, and its value in the test counterpart is 
$\mathcal L_{\rm test} = 3.98\,\%$. As an illustration, three examples of 
train and test solutions are depicted in 
\figurename~\ref{fig:rrae-sol-reconstruction}.

\subsubsection{Direct problem}

Regarding the methodology for this 1d field QoI, the next step consists in 
finding a mapping relating geometries latent coefficients $\bm A$ and 
solutions latent coefficients $\bm B$. Again, an MLP is trained to this end, 
and the sought mapping can be expressed as 
$\mathcal{NN}_{\rm 1d}^d:\bm \alpha_j\rightarrow \bm \beta_j$, where the 
superscript $d$ stands for direct problem. Recall that $\bm \alpha_j\in \mathbb R^2$ 
and $\bm \beta_j\in \mathbb R^1$ due to the chosen $k_{\max}$ 
(hyper-)parameters in previous sections.

Analogously, the hyper-parameters of the MLP architecture are displayed in 
\tablename~\ref{tab:nnd-1d}, the optimization parameters are highlighted in 
\tablename~\ref{tab:nnd-1d-optimization-parameters}, and both input and output 
are normalized before performing the training. The loss value in the train is 
$\textrm{MSE}_{\rm train} = 2.2 \cdot 10^{-3}$. The coefficient of 
determination $R^2$ is another metric to assess the goodness of the fit, 
having obtained $R^2_{\rm train} = 0.998$ and $R^2_{\rm test} = 0.989$. The 
prediction of every latent coefficient $\beta_0$ through 
$\mathcal{NN}_{\rm 1d}^d$ compared with its ground-truth value is depicted in 
\figurename~\ref{fig:nnd-1d-latent}. Note the smoothness of latent 
coefficient, achieved thanks to the use of MLP-RRAE of the solutions. In other 
words, the dimensionality reduction keeps the physical meaning of varying the 
volume fraction $f$, making the relation between latent spaces an easier task 
with respect to the previous scalar counterpart.

Now, we have three trained models: a CNN-RRAE for the geometries, an MLP-RRAE 
for 1d-solutions, and an MLP $\mathcal{NN}_{\rm 1d}^d$ relating latent 
coefficients of the both. Joining all them properly by following the procedure 
outlined in Section~\ref{sec:methodology}, we can generate the pipeline 
$\texttt{geo2sol}_{\rm 1d}$, whose input is a geometry $\bm x_j$, and its 
output is $\tilde{\bm s}_j$, the von Mises stress distribution across the  
diagonal (top-left to bottom-right elements) of such geometry i.e., the 1d 
field highlighted in this section. To assess the accuracy of predictions, the 
coefficient of determination $R^2$ between ground-truth and predicted 1d 
curves is computed per sample, and then averaged across the dataset. These 
coefficients are:
\begin{equation*}
R^2_{\rm train} = 0.997, \quad R^2_{\rm test} = 0.997.
\end{equation*}

\begin{figure}[H]
    \centering
    \subfloat[Train example]{\includegraphics[width=0.6\linewidth]{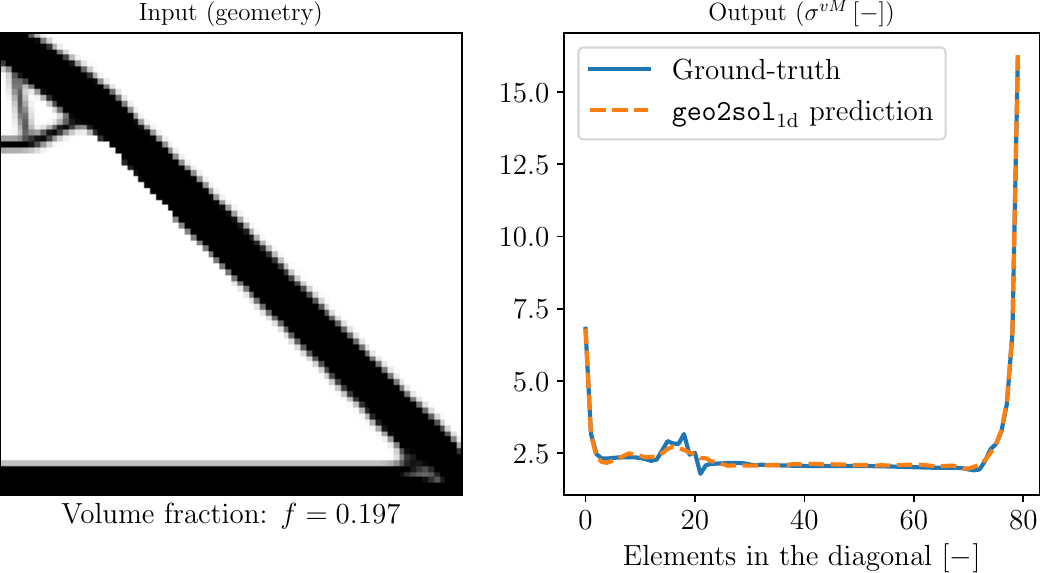}}\\
    \subfloat[Test example]{\includegraphics[width=0.6\linewidth]{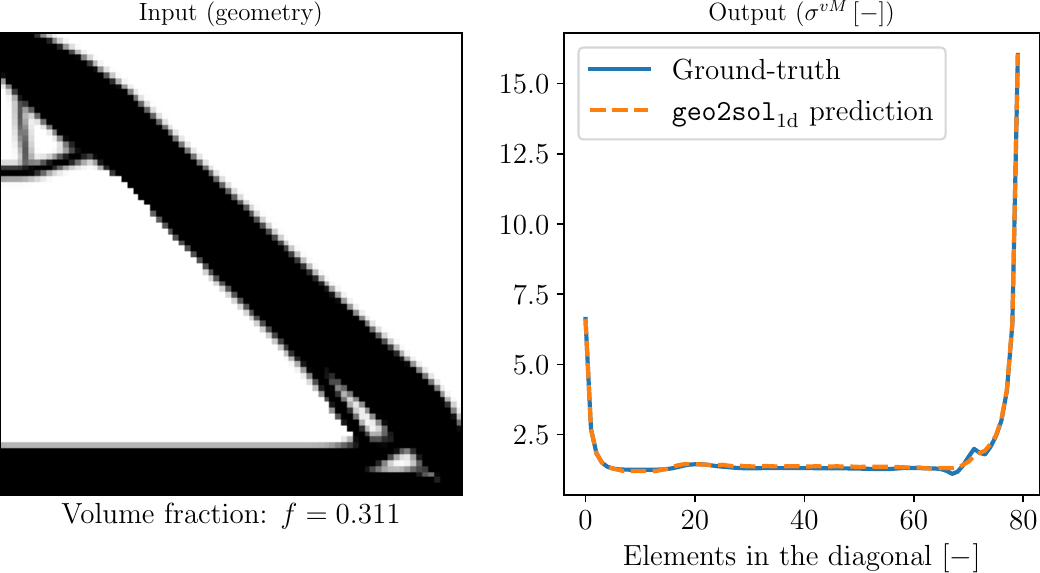}}
    \caption{$\texttt{geo2sol}_{\rm 1d}$ example predictions in (a) train and (b) test.}
    \label{fig:geo2sol1d}
\end{figure}

\figurename~\ref{fig:geo2sol1d} displays samples of successful $\texttt{geo2sol}_{\rm 1d}$ predictions in train and test sets.

\subsubsection{Inverse problem}

This procedure is analogous to the previous subsection, but generating an 
inverse mapping relating latent coefficients i.e., from $\bm B$ to $\bm A$. 
The mapping $\mathcal{NN}^i_{\rm 1d}$, where now $i$ stands for inverse, is 
again an MLP that inputs $\beta_j\in \mathbb R^1$ and outputs 
$\bm \alpha_j \in\mathbb R^2$. Their hyper-parameters are the same as the 
previous (\tablename~\ref{tab:nnd-1d}), with the input and output dimensions 
permuted, but with the same hidden layer configurations. The same training 
parameters (\tablename~\ref{tab:nnd-1d-optimization-parameters}) are likewise 
used. The loss value is ${\rm MSE}_{\rm train} = 1.1\cdot 10^{-3}$, and the 
coefficients of determination, $R^2_{\rm train} = 0.999$ and 
$R^2_{\rm test} = 0.999$. The two curves for the prediction of coefficients 
$\bm\alpha = [\alpha_0, \alpha_1]$ are depicted in 
\figurename~\ref{fig:nni-1d-latent}.

With $\mathcal{NN}_{\rm 1d}^i$, the $\texttt{sol2geo}_{\rm 1d}$ pipeline can 
be built, whose input is the von Mises stress curve (across the diagonal of 
the domain $\Omega$), predicting the geometry which generates such 
curve---hence the inverse problem is addressed. The encoder of MLP-RRAE of 
solutions maps the input curve $\bm s_j$ into a (single) latent coefficient 
$\beta_{0,j}$. Then, such coefficient is mapped through 
$\mathcal{NN}_{\rm 1d}^i$ to generate the latent geometry coefficients 
$\bm \alpha_j$, to be further mapped through the CNN-RRAE decoder, thus 
generating a geometry $\tilde{\bm x}_j$ as the output. Evaluating the 
reconstruction loss, Equation~\eqref{eq:reconstruction-loss},  between 
ground-truth an $\texttt{sol2geo}_{\rm 1d}$ predicted geometries yields the 
following:
\begin{equation*}
\mathcal{L}_{\rm train} = 6.38\,\%, \quad \mathcal{L}_{\rm test} = 7.63\,\%.
\end{equation*}

\begin{figure}[H]
    \centering
    \subfloat[Train example]{\includegraphics[width=.8\linewidth]{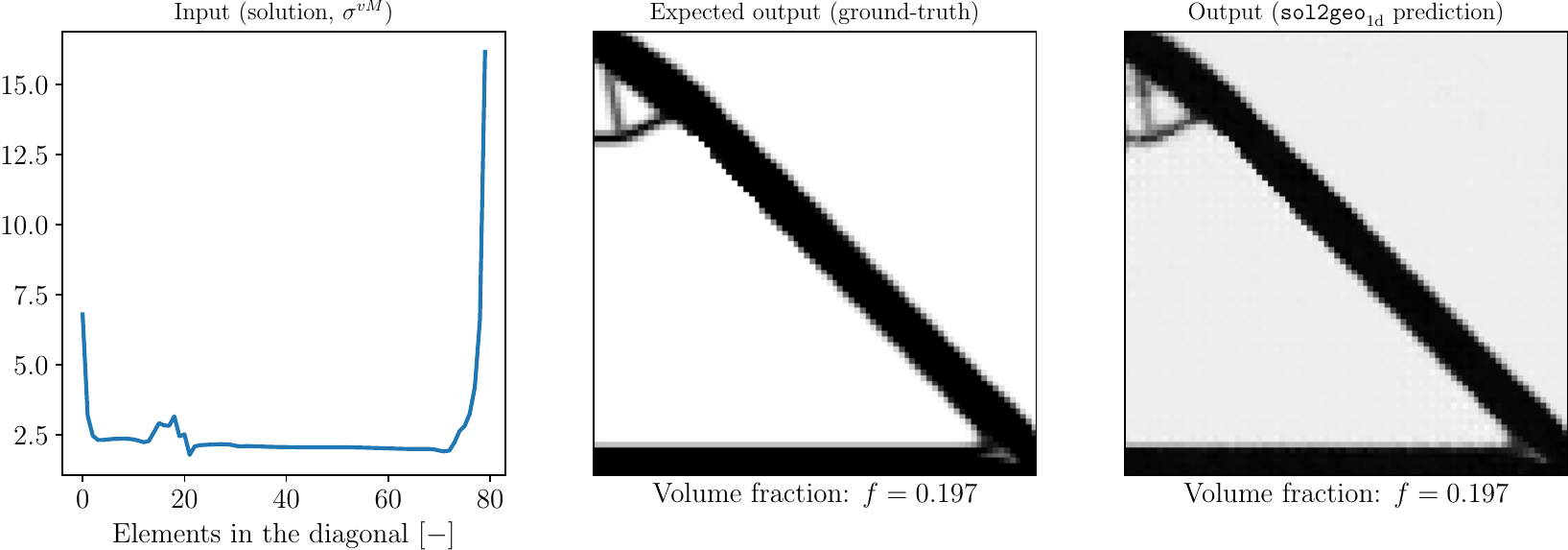}}\\
    \subfloat[Test example]{\includegraphics[width=.8\linewidth]{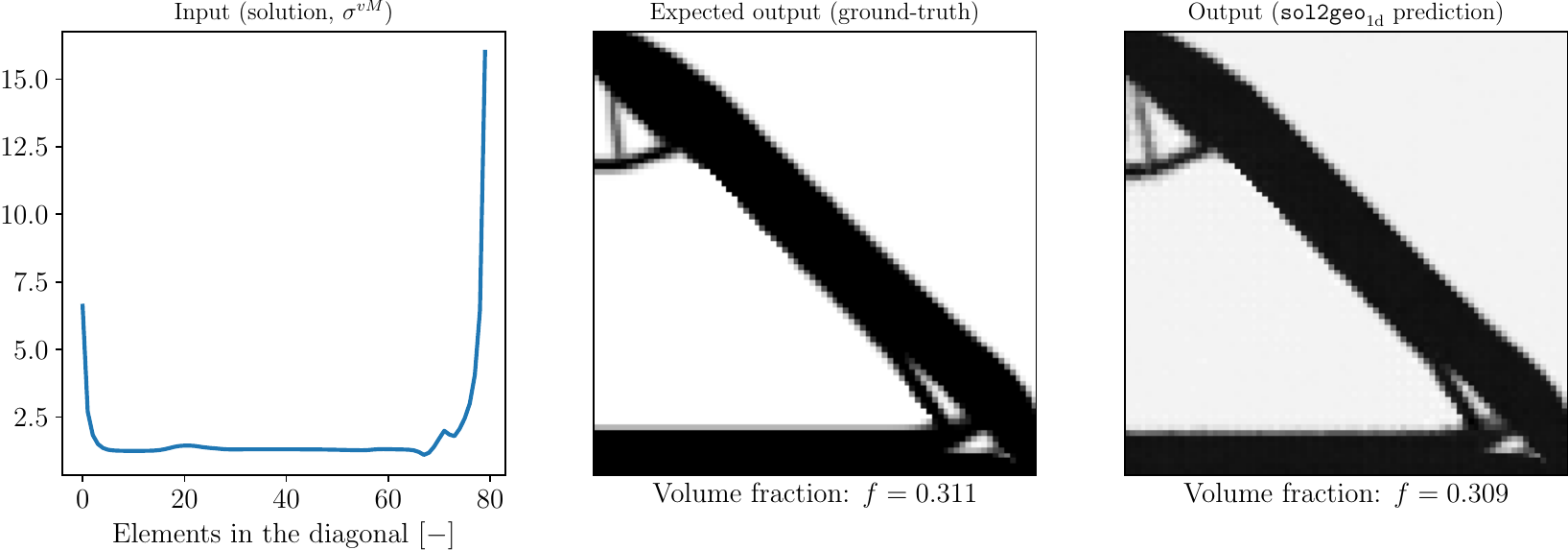}}
    \caption{$\texttt{sol2geo}_{\rm 1d}$ example predictions in (a) train and (b) test. Note the closeness of volume fractions in ground-truth and predicted geometries -- exact in the train -- despite not having prescribed such constraint in the pipeline.}
    \label{fig:sol2geo1d}
\end{figure}

\figurename~\ref{fig:sol2geo1d} displays samples of successful $\texttt{sol2geo}_{\rm 1d}$ predictions in train and test sets.
\begin{figure}[H]
    \centering
    \includegraphics[width=0.6\linewidth]{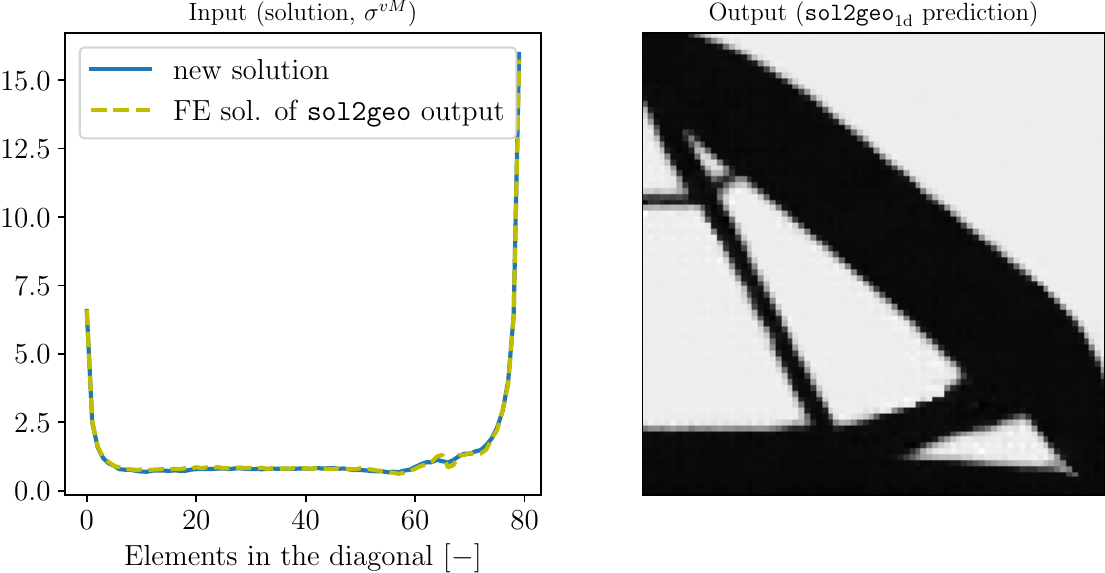}
    \caption{Example to assess the validity of solutions from generated geometries. A new solution (blue, left) is fed to the $\texttt{sol2geo}_{\rm 1d}$ pipeline, yielding a geometry (right). Then, this geometry is solved via FE, and the $\sigma^{vM}$ of the diagonal elements i.e., the high-fidelity solution, is represented against the input (dashed yellow, left).}
    \label{fig:sol2geo1d-interp}
\end{figure}

Lastly in this section, when predicting geometries, one can ask whether these 
geometries are physically reasonable. This is what we intend to evaluate to 
conclude this section. To do so, a new solution is going to be tested by 
comparing it with the actual solution of the (new) generated geometry. To be 
new, the solution is generated as an interpolation in the latent space of 
solutions, therefore, the geometry produced in the $\texttt{sol2geo}_{\rm 1d}$ 
pipeline will be new---never seen before by the model. This geometry is solved 
via FE, so a high-fidelity curve related to this geometry is obtained. 
Finally, this is compared with the input, i.e. the new solution (result of 
interpolation). This idea is displayed in 
\figurename~\ref{fig:sol2geo1d-interp}.

\subsection{2d field quantity of interest in a fixed domain}

The last case in this paper showcases the direct and inverse problems when the solutions (QoI) are a 2d field defined in the fixed domain $\Omega$, that is:
\begin{equation*}
\bm s_j := \bm \sigma^{vM,*}_j,
\end{equation*}
properly rearranged in a matrix i.e., replicating the 2d mesh. In this case, full resolution $\bm \sigma^{vM}$ solutions are predicted using this methodology.

\begin{figure}[H]
    \subfloat[Train]
    {\includegraphics[width=.45\linewidth]{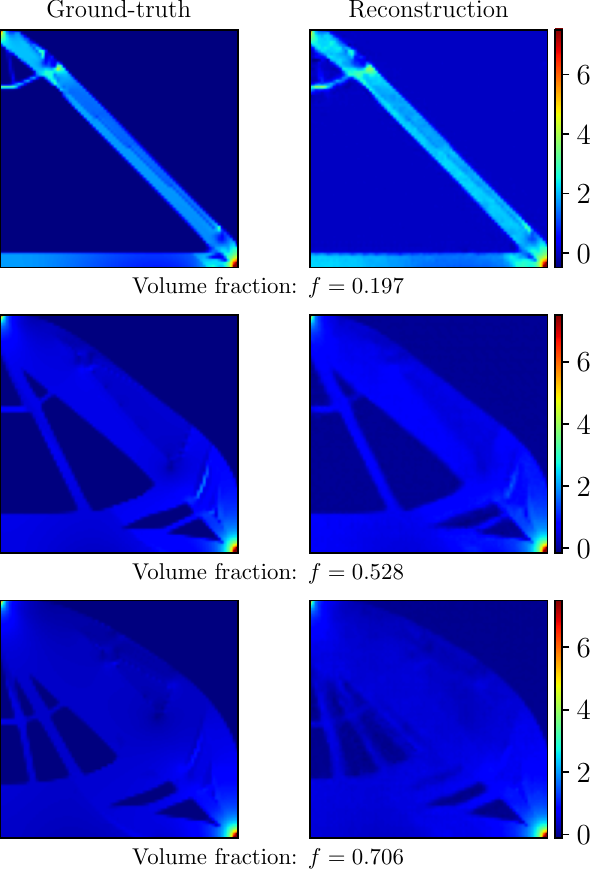}}\hfill%
    \subfloat[Test]
    {\includegraphics[width=.45\linewidth]{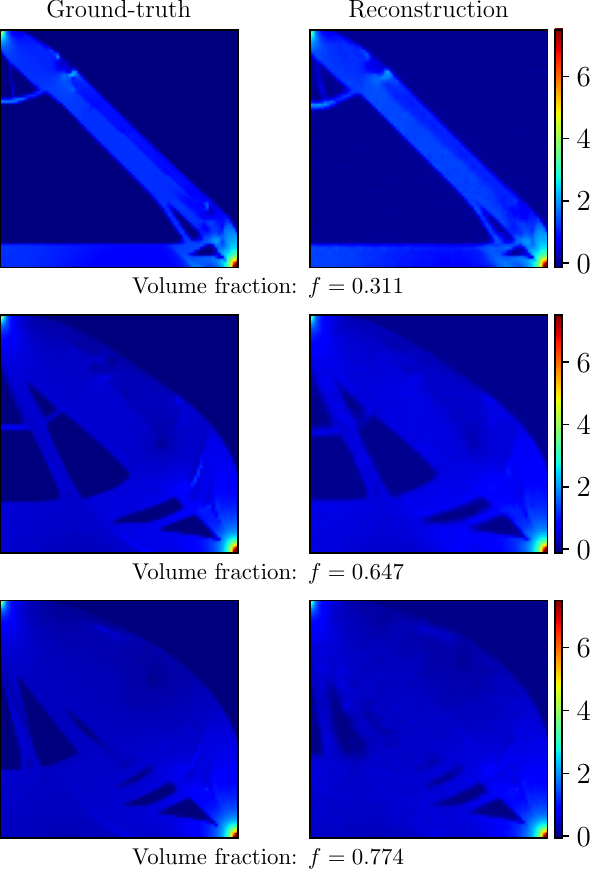}}
    \caption{CNN-RRAE 2d solutions (von Mises stress distribution) model. Example reconstructions in (a) train and (b) test sets. For the sake of visualization, the maximum value in the plot is $7.5$.}
    \label{fig:rrae-sol2d-reconstruction}
\end{figure}

On the one hand, the already trained CNN-RRAE highlighted in 
Section~\ref{subsec:scalar-qoi} is used as the geometries model. On the other 
hand, a new CNN-RRAE for solutions is built and trained. The architecture is 
the same (\tablename~\ref{tab:geometries-rrae}), but using a $k_{\max} = 1$ in 
order to compress in a unique latent coefficient. Training the model with 
\tablename~\ref{tab:geometries-rrae-optimization-parameters} parameters, the 
reconstruction loss of the solutions RRAE is 
$\mathcal L_{\rm train} = 12.23\,\%$ in the training set, and its value in the 
test counterpart is $\mathcal L_{\rm test} = 13.87\,\%$. Examples for both 
train and test reconstructions are displayed in 
\figurename~\ref{fig:rrae-sol2d-reconstruction}.

\subsubsection{Direct problem}

The direct problem is addressed by building a mapping $\mathcal{NN}_{\rm 2d}^d$ 
relating the geometry latent coefficients $\bm A$ and the new 2d-field 
solution latent coefficients $\bm B$---the same symbol is used, with a slight 
abuse of notation. Since the architecture is preserved for all the mappings, 
the same hyper-parameters are used (\tablename~\ref{tab:nnd-1d} and 
\tablename~\ref{tab:nnd-1d-optimization-parameters}).

The loss value is $\textrm{MSE}_{\rm train} = 1.5\cdot 10^{-4}$, and 
coefficients of determination are $R^2_{\rm train} = 1.000$, 
$R^2_{\rm test} = 0.999$. Note that since the solution part has more 
information in this 2d case i.e., it is not compressed like the 1d or scalar 
case, it is easier to relate the latent coefficients to the solution part, and 
hence the metrics are better. The latent coefficient $\beta_0$ prediction 
results of $\mathcal{NN}_{\rm 2d}^d$ are depicted in 
\figurename~\ref{fig:nnd-2d-latent}.

In this subsection, the pipeline $\texttt{geo2sol}_{\rm 2d}$ is built: 
CNN-RRAE geometries model encoding, map the latent coefficients through 
$\mathcal{NN}_{\rm 2d}^d$, and decode through the CNN-RRAE 2d solutions model. 
Once built, the evaluation of the reconstruction loss, 
Equation~\eqref{eq:reconstruction-loss},  between ground-truth an 
$\texttt{geo2sol}_{\rm 2d}$ predicted geometries yields the following:
\begin{equation*}
\mathcal{L}_{\rm train} = 10.39\,\%, \quad \mathcal{L}_{\rm test} = 11.59\,\%.
\end{equation*}

\begin{figure}[H]
    \centering
    \subfloat[Train example.]
    {\includegraphics[width=.8\linewidth]{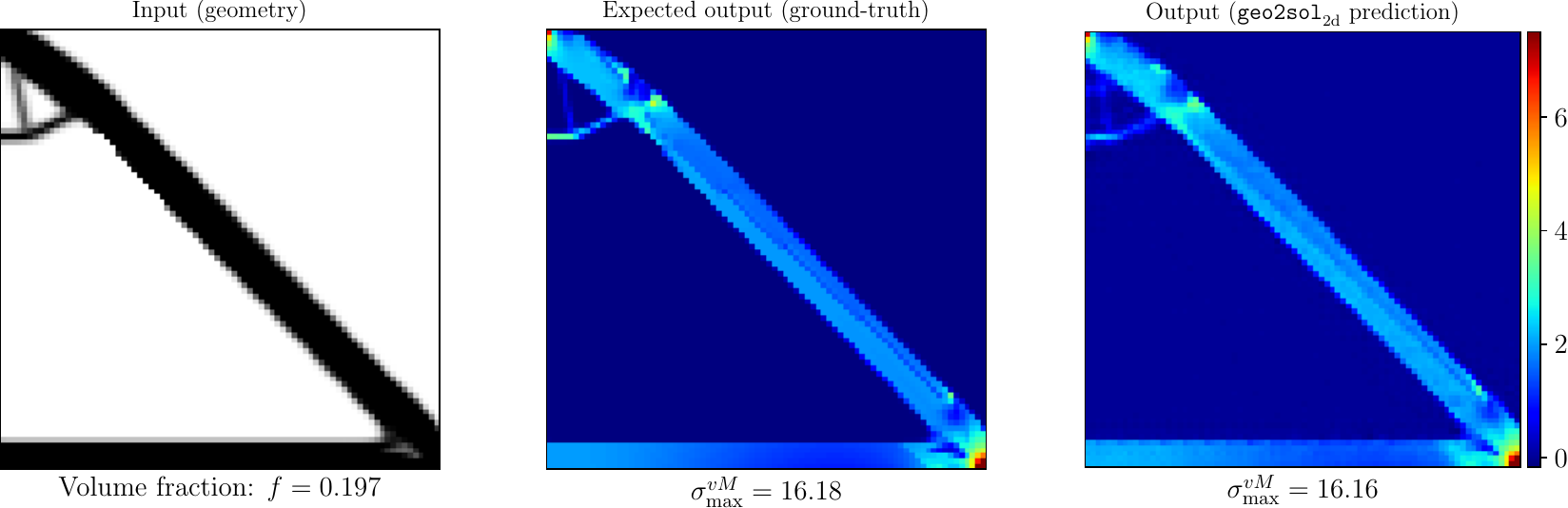}}\\
    \subfloat[Test example.]
    {\includegraphics[width=.8\linewidth]{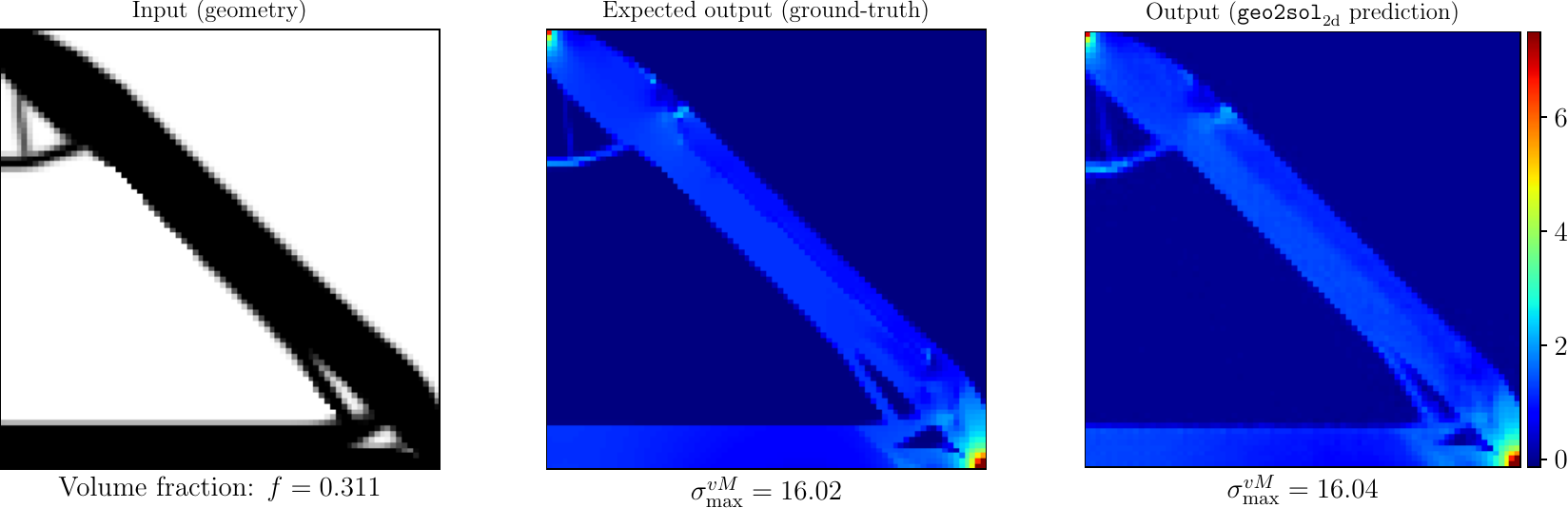}}
    \caption{$\texttt{geo2sol}_{\rm 2d}$ example predictions in (a) train and (b) test. Note how close are the maximum stress values in ground-truth and predicted stress fields. The maximum value in the plot is 7.5 again.}
    \label{fig:geo2sol2d}
\end{figure}

Two examples (train, test) are shown in \figurename~\ref{fig:geo2sol2d}. Note 
the accuracy of the prediction in how close the maximum stress values 
$\sigma^{vM}_{\max}$ are in both ground-truth and $\texttt{geo2sol}_{\rm 2d}$ 
predictions. This may serve as an alternative to the direct scalar QoI case 
($\texttt{geo2sol}_s$).

\subsubsection{Inverse problem}

The last mapping $\mathcal{NN}_{\rm 2d}^i$ is built to map $\bm B$ 2d solution 
latent coefficients into $\bm A$ geometry latent coefficients. Using the same 
MLP architecture and training parameters as the previous cases, the loss value 
obtained is $\mathrm{MSE}_{\rm train}=6.3\cdot 10^{-4}$, and the coefficients 
of determination $R^2_{\rm train} = 1.000$, $R^2_{\rm test} = 0.999$. The 
fitting curves are depicted in \figurename~\ref{fig:nni-2d-latent}.

Now, the inverse mapping $\texttt{sol2geo}_{\rm 2d}$ with respect to the 
previous subsection is built. This implies taking the stress distribution in 
the mesh as input $\bm s_j$, encode it via the CNN-RRAE 2d solutions model, 
map the corresponding $\beta_{0,j}$ coefficient through $\mathcal{NN}_{\rm 2d}^i$ 
into $\bm \alpha_j$, and decode it with the CNN-RRAE geometries model to 
obtain the predicted geometry $\tilde{\bm x}_j$. Evaluating the reconstruction 
loss in the pipeline, between ground-truth an $\texttt{sol2geo}_{\rm 2d}$ 
predicted geometries, yields the following:
\begin{equation*}
\mathcal{L}_{\rm train} = 5.74\,\%, \quad \mathcal{L}_{\rm test} = 7.80\,\%.
\end{equation*}

\begin{figure}[H]
    \centering
    \subfloat[Train example.]
    {\includegraphics[width=.8\linewidth]{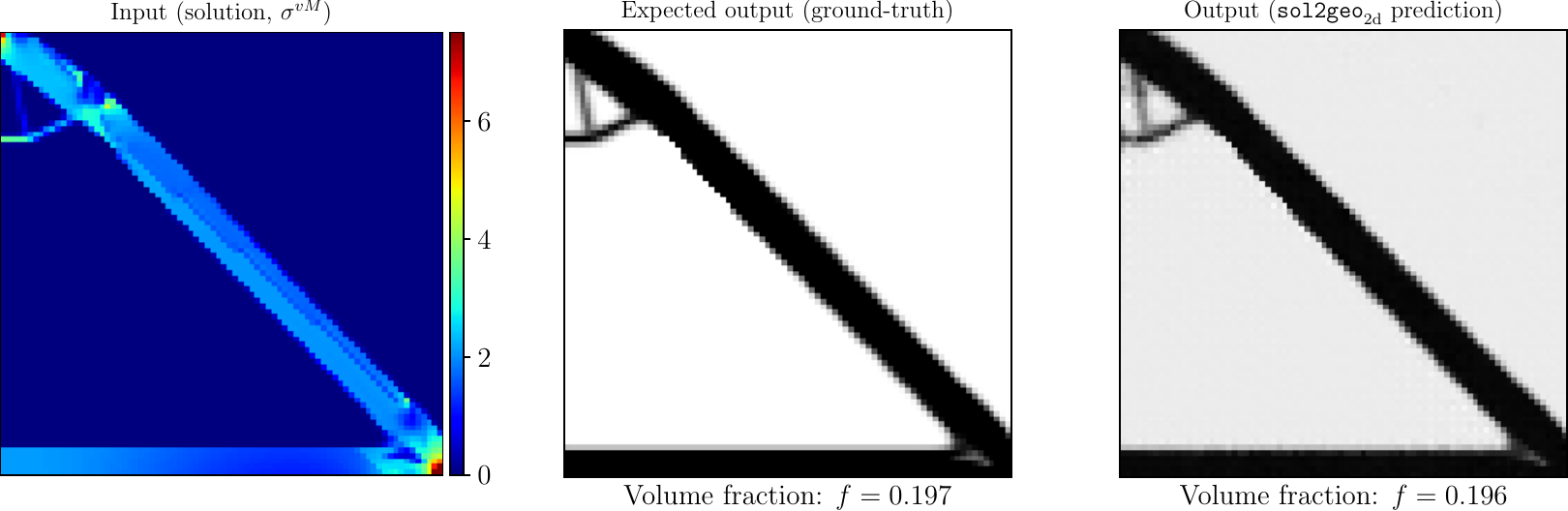}}\\
    \subfloat[Test example.]
    {\includegraphics[width=.8\linewidth]{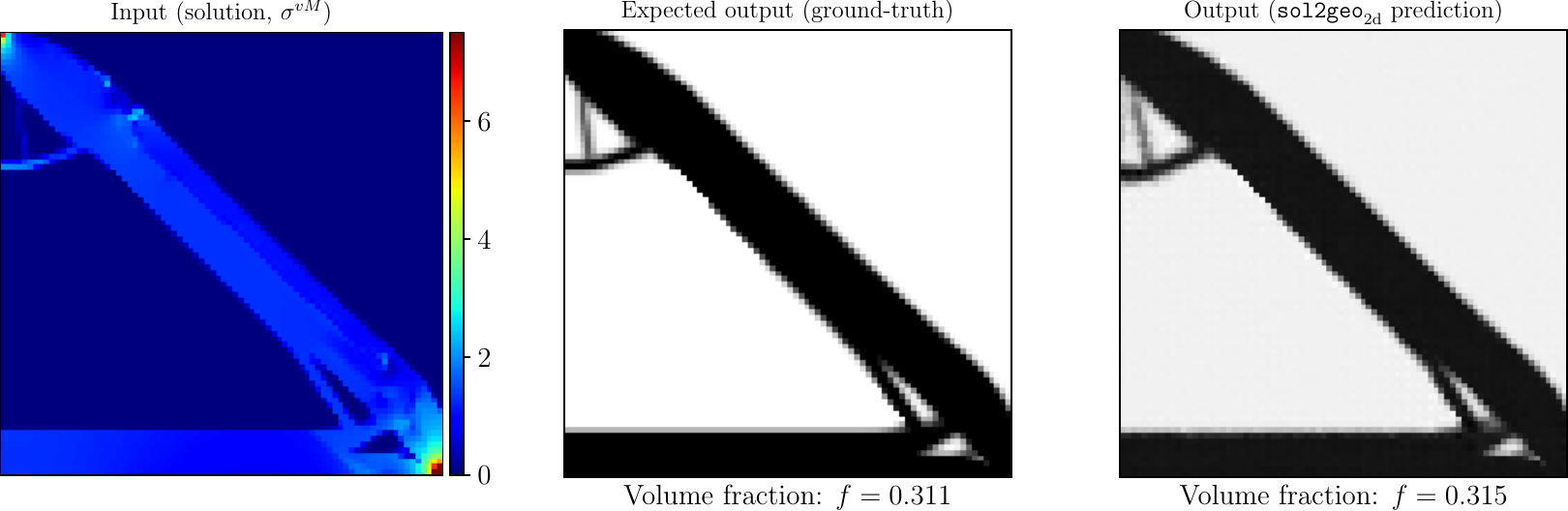}}
    \caption{$\texttt{sol2geo}_{\rm 2d}$ example predictions in (a) train and (b) test. Note how close are the volume fractions in ground-truth and predicted geometries. The maximum value in the plot is 7.5 again.}
    \label{fig:sol2geo2d}
\end{figure}

Examples of train and test showcasing this pipeline are depicted in \figurename~\ref{fig:sol2geo2d}.

\section{Conclusions and Perspectives}\label{sec:concl_and_persp}

This work has introduced a data-driven framework that combines Rank Reduction 
Autoencoders (RRAEs) with neural latent-space mappings to address both forward 
and inverse problems in topology optimization (TO), both of them predicted via regression.
By considering a low-rank approximation of the latent space, the proposed methodology enables efficient 
dimensionality reduction while preserving the dominant physical trends 
governing optimized geometries and their mechanical responses.

The numerical results demonstrate that the effectiveness of the framework 
strongly depends on the richness of the selected quantity of interest (QoI). 
When the response is compressed into a scalar metric, such as the maximum von 
Mises stress, the resulting geometry-solution mapping becomes inherently 
non-bijective or ill-conditioned. This limitation is reflected in reduced 
predictive accuracy, increased sensitivity to extrapolation, and overfitting 
tendencies in the inverse problem. These observations highlight the intrinsic 
difficulty of recovering detailed geometric information from highly localized 
data.
In contrast, when the QoI is defined as a one-dimensional or two-dimensional 
stress field, the proposed approach achieves a substantial improvement in both 
forward and inverse predictions. The use of an RRAE for the solution fields 
yields smooth, physically meaningful latent coefficients that correlate well 
with variations in the design parameters. As a result, the learned mappings 
between geometry and solution latent spaces exhibit high accuracy and 
robustness, enabling reliable reconstruction of geometries and responses 
without enforcing explicit constraints such as volume fraction preservation.
The best overall performance is obtained when the full two-dimensional von 
Mises stress field is considered as the QoI. In this case, the increased 
informational content of the solution leads to highly accurate latent-space 
regressions and consistent predictions in both direct and inverse settings. 
Despite the higher complexity of the data, the RRAE-based compression remains 
effective, and the reconstructed stress fields and geometries closely match 
their high-fidelity counterparts. These results suggest that the proposed 
framework can serve as an efficient surrogate for Finite Element analyses in 
TO workflows.

Beyond surrogate modeling, the methodology naturally supports Generative 
Design (GD) capabilities through latent-space interpolation due to the use of the 
SVD model-order reduction, enabling the synthesis of new geometries and 
responses that remain physically plausible. Although this aspect was only 
preliminarily explored, it highlights the potential of RRAEs as a unifying 
tool for reduced-order modeling, inverse design, and data-driven exploration or 
GD of optimized structures.

Future work will focus on extending the framework to multi-parameter design 
spaces, incorporating varying boundary conditions and loading scenarios, and 
enforcing physics-informed constraints within the latent mappings. 
Additionally, integrating the proposed approach into optimization loops and 
uncertainty-aware design pipelines represents a promising direction for 
advancing data-driven generative mechanical design.

\appendix \setcounter{table}{0}\setcounter{figure}{0}
\section{MLP models details}

\subsection[Architecture and optimization hyper-parameters]{$\mathcal{NN}_{\rm 1d}^d$ and $\mathcal{NN}_{\rm 1d}^i$ architecture and optimization hyper-parameters}\label{app:nnd-1d}

\begin{table}[H]
    \centering
    \begin{tabular}{c|cccc}
    \hline
    Module & Layer & Input shape & Output shape & Activation \\ \hline
    Input layer & Dense & $(n_b, 2)$ & $(n_b, 8)$ & ELU \\ 
    \rowcolor{lightgray}Hidden layer \#1 & Dense & $(n_b, 8)$ & $(n_b, 16)$ & ELU \\
    Hidden layer \#2 & Dense & $(n_b, 16)$ & $(n_b, 8)$ & ELU \\
    \rowcolor{lightgray}Hidden layer \#3 & Dense & $(n_b, 8)$ & $(n_b, 4)$ & ELU \\
    Output layer & Dense & $(n_b, 4)$ & $(n_b, 1)$ & Linear \\ 
    \hline
    \end{tabular}
    \caption{Architecture of $\mathcal{NN}_{\rm 1d}^d$ (and 
    $\mathcal{NN}_{\rm 1d}^i$, permuting input and output shapes), the MLP 
    model relating latent coefficients from geometry to solution in the 1d QoI 
    case. }
    \label{tab:nnd-1d}
\end{table}

\begin{table}[H]
    \centering
    \begin{tabular}{c|c}
         Parameter & Choice \\ \hline
         Optimizer & nadam \\
         \rowcolor{lightgray}Learning rate & $10^{-3}$\\
         Epochs & $3000$ \\
         \rowcolor{lightgray}Batch size & $100$ \\
         Loss & MSE 
    \end{tabular}
    \caption{$\mathcal{NN}_{\rm 1d}^d$ and $\mathcal{NN}_{\rm 1d}^i$ MLP models. Optimization parameters of the training.}
    \label{tab:nnd-1d-optimization-parameters}
\end{table}

\subsection[Fitting curves]{$\mathcal{NN}_{\star}^d$ and $\mathcal{NN}_{\star}^i$ fitting curves}\label{app:mlp-fitting-curves}

\begin{figure}[H]
    \centering
    \includegraphics[width=0.5\linewidth]{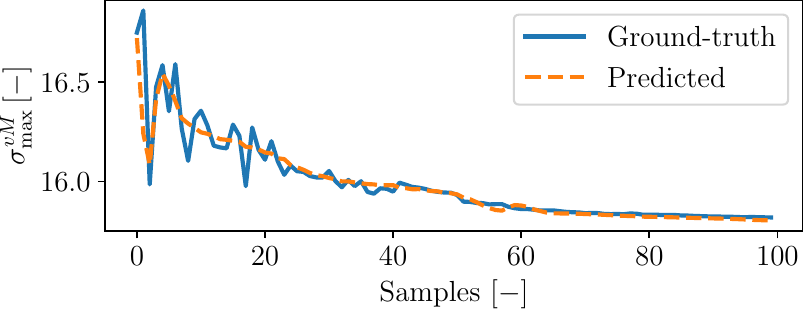}
    \caption{$\mathcal{NN}_{s}^d$ MLP model. Prediction of the scalar QoI $\sigma_{\max}^{vM}$.}
    \label{fig:nnd-s-latent}
\end{figure}

\begin{figure}[H]
    \centering
    \includegraphics[width=\linewidth]{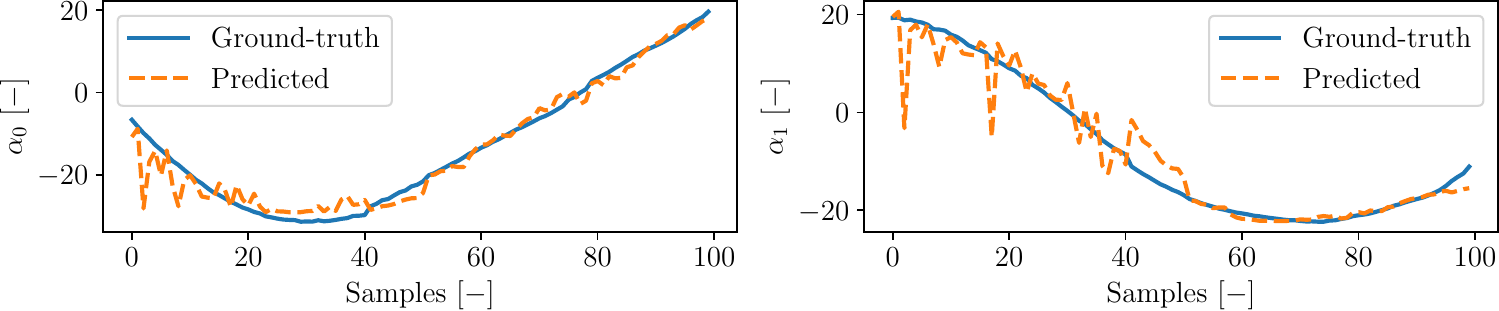}
    \caption{$\mathcal{NN}_{s}^i$ MLP model. Prediction of the geometries latent coefficients $\alpha_0$ and $\alpha_1$.}
    \label{fig:nni-s-latent}
\end{figure}

\begin{figure}[H]
    \centering
    \includegraphics[width=0.5\linewidth]{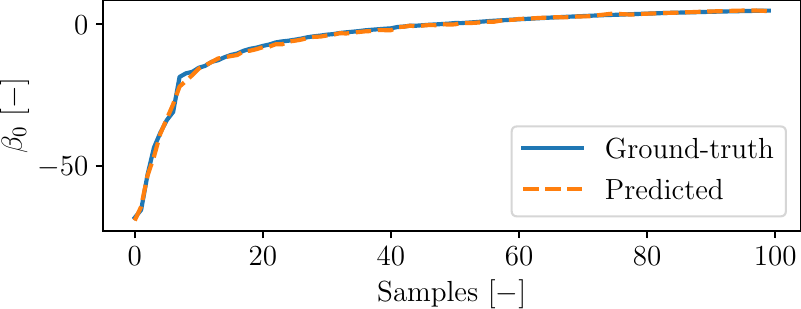}
    \caption{$\mathcal{NN}_{\rm 1d}^d$ MLP model. Prediction of the solutions latent coefficients $\beta_0$.}
    \label{fig:nnd-1d-latent}
\end{figure}

\begin{figure}[H]
    \centering
    \includegraphics[width=\linewidth]{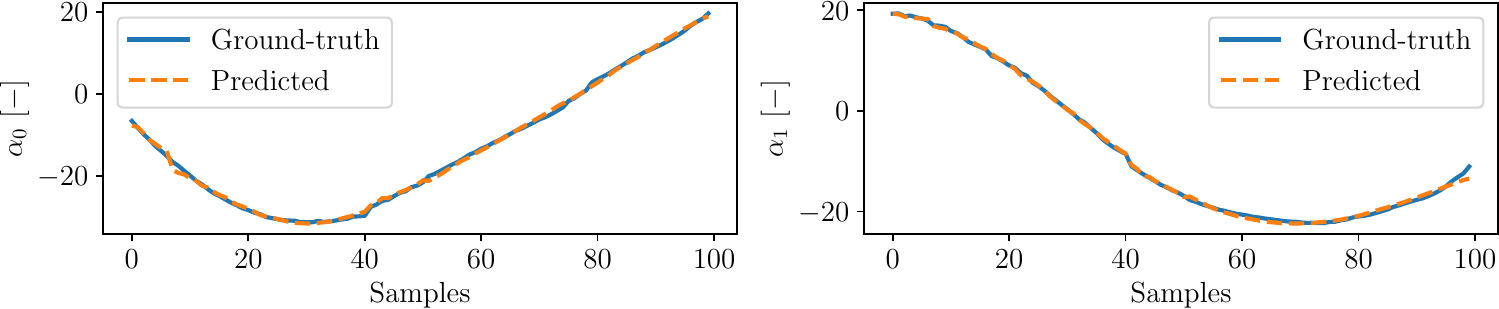}
    \caption{$\mathcal{NN}_{\rm 1d}^i$ MLP model. Prediction of the geometries latent coefficients $\alpha_0$ and $\alpha_1$.}
    \label{fig:nni-1d-latent}
\end{figure}

\begin{figure}[H]
    \centering
    \includegraphics[width=0.5\linewidth]{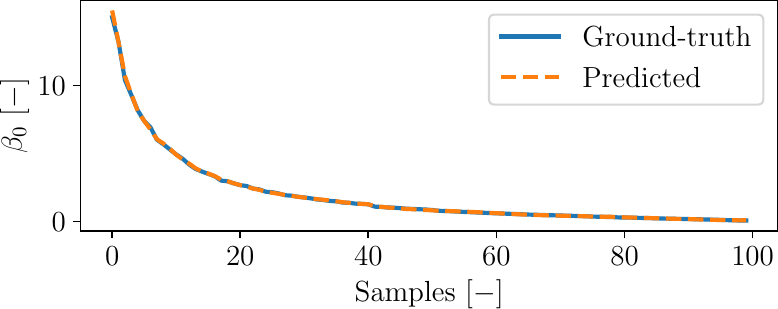}
    \caption{$\mathcal{NN}_{\rm 2d}^d$ MLP model. Prediction of the solutions latent coefficients $\beta_0$.}
    \label{fig:nnd-2d-latent}
\end{figure}

\begin{figure}[H]
    \centering
    \includegraphics[width=\linewidth]{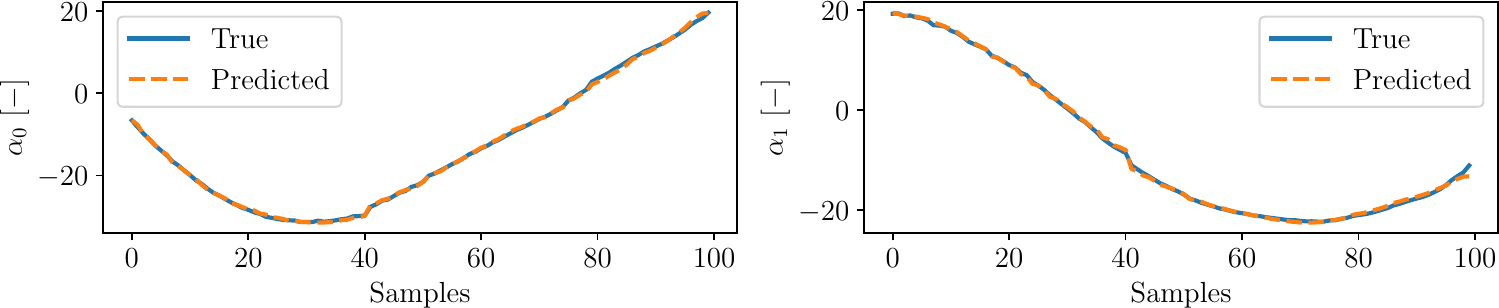}
    \caption{$\mathcal{NN}_{\rm 2d}^i$ MLP model. Prediction of the geometries latent coefficients $\alpha_0$ and $\alpha_1$.}
    \label{fig:nni-2d-latent}
\end{figure}

\bibliographystyle{unsrt}
\bibliography{rraes-topology-optimization}

\end{document}